\documentclass[11pt,a4paper]{article}
\usepackage{epsfig}
\usepackage{amssymb}
\usepackage{latexsym}
\usepackage{epic}
\usepackage{eepic}
\usepackage{amsmath}
\usepackage{epstopdf}

\newcommand{\ep}{\epsilon}
\newcommand{\bep}{b_{\ep}(d)}
\newcommand{\Ex}{\mathrm{E}}
\newcommand{\Var}{\mathrm{Var}}
\newcommand{\Cov}{\mathrm{Cov}}
\newcommand{\qt}{\tilde{q}}
\newcommand{\pt}{\tilde{p}}
\newcommand{\Qt}{\tilde{Q}}

\newcommand{\mc}{\mathcal}

\newcommand{\prob}{\stackrel{\rm P}{\to}}
\newcommand{\Or}{{\rm O}}

\newcommand{\oor}{{\rm o}}
\newcommand{\ninf}{\mbox{ as } n \to \infty}
\newtheorem{theorem}{Theorem}
\newtheorem{lemma}{Lemma}
\newtheorem{proposition}{Proposition}

\newenvironment{reflist}{\begin{list}{}
         {\itemsep=20pt \parsep=3pt
          \topsep=0pt  \parskip=20pt  \listparindent=-.15in
           \leftmargin= 0.15in }
         \item \ \vspace{-.35in} }
         {\end{list}}

\begin{document}

\author{David K\"{a}llberg$^{1}$, Nikolaj Leonenko$^{2}$, Oleg\ Seleznjev$%
^{1}$ \vspace{0.5cm}\hfill \\
$^1 \,$Department of Mathematics and Mathematical Statistics\\
Ume{\aa } University, SE-901 87 Ume\aa , Sweden,\\
$^2 \,$ School of Mathematics, Cardiff University,\\
Senghennydd Road, Cardiff CF24 4YH, UK\\
}
\date{  }
\title{Statistical estimation of quadratic R\'{e}nyi entropy for a
stationary $m$-dependent sequence }
\maketitle

\begin{abstract}
The R\'{e}nyi entropy  is a generalization of the Shannon entropy  and is widely used in mathematical statistics and applied sciences
for quantifying the uncertainty in  a probability
distribution.
We consider estimation of the quadratic R\'{e}nyi entropy and related functionals
for the marginal distribution of a stationary $m$-dependent sequence.
The $U$-statistic estimators under study are based on the number of $\ep$-close vector observations in the corresponding  sample.
A variety of asymptotic properties for these estimators are obtained (e.g., consistency, asymptotic normality,
Poisson convergence). The results can be used in diverse statistical and computer science problems
whenever the conventional independence assumption is too strong (e.g., $\ep$-keys
 in time series databases, distribution
identification problems for dependent samples).
\end{abstract}

 \baselineskip=3.4 ex
\noindent \emph{AMS 2010 subject classification:} 62G05, 62G20, 62M99, 94A17
\medskip

\noindent \emph{Keywords:}
Entropy estimation, quadratic R\'{e}nyi entropy, stationary $m$-dependent sequence, inter-point distances, $U$-statistics

\section{ Introduction}

Entropy is applied  in information theory and statistics
for characterizing the diversity or uncertainty in a probability distribution.
For a continuous distribution $\mc{P}$ with density $p(x),x\in R^{d}$, the R\'{e}nyi entropy is defined (R\'{e}nyi, 1970)
as
\begin{equation*}
h_{s }(\mc{P}):=\frac{1}{1-s }\log \left(
\int_{R^{d}}p(x)^{s }\,dx\right) ,\quad s \neq 1.
\end{equation*}%
Henceforth we use $\log x$ to denote the natural logarithm of $x$. The R\'{e}nyi entropy $h_s$ is a generalization of the Shannon entropy (Shannon,
1948):
\begin{equation*}
h_{1}(\mc{P})=\lim_{s \rightarrow 1}h_{s }(\mc{P})=-\int_{R^{d}}p(x)\log p(x)\,dx.
\end{equation*}
From the statistical point of view, the quadratic R\'{e}nyi entropy
\begin{equation*}
h_{2}=h_{2}(\mc{P})=-\log \left( \int_{R^{d}}p(x)^{2}\,dx\right) ,
\end{equation*}
is the simplest point on the R\'{e}nyi spectrum $\{h_{s }(\mc{P}
),s \in \mc S \}$, where $\mc S $ is subset of $R=R^{1}$, such that
entropies exist. Note that
$h_{2}=-\log (q_{2})$,
where we assume that the quadratic functional
\begin{equation*}
q_{2}=q_{2}(\mc{P}):=\int_{R^{d}}p(x)^{2}\,dx.
\end{equation*}
is well defined, and hence the point $s =2$ belongs to the set $\mc S$.
More entropy generalizations are known in information theory, e.g., the Tsallis entropy  (Tsallis, 1988):
\begin{equation*}
T_{s}(\mc{P}):=\frac{1}{1-s}(1-\int_{R^{d}} p(x)^{s}dx),\quad s\neq 1.
\end{equation*}
The R\'{e}nyi entropy (or information) for stationary processes can be
understood as that of the corresponding ergodic or marginal distributions, see,
e.g.,  Gregorio and Iacus (2003), where the R\'{e}nyi
entropy is computed for a large class of ergodic diffusion
processes.

Numerous  applications of the R\'{e}nyi entropy in information theoretic learning,
statistics (e.g., classification, distribution identification problems, statistical inference), computer science
(e.g., average case analysis for random databases, pattern recognition, image matching), and econometrics
are discussed, e.g., in Principe (2010), Kapur (1989), Kapur and Kesavan
(1992), Pardo (2006), Escolano et al.\ (2009),  Neemuchwala et al.\ (2005),  Ullah
(1996),
Baryshnikov et al.\ (2009), Seleznjev and Thalheim (2003, 2010), Thalheim
(2000), Leonenko et al.\ (2008), and Leonenko and Seleznjev (2010).

Various estimators for the quadratic functional $q_2$ and the entropy $h_2$
for {\it independent} samples have been studied.
Leonenko et al.\ (2008) obtain consistency of nearest-neighbor estimators for $h_{2}$,
see also Penrose and Yukich (2011) and the references therein.
Bickel and Ritov (1988) and Gin\'{e} and Nickl (2008) show rate optimality, efficiency, and asymptotic normality of kernel-based estimators for $q_2$ in the one-dimensional case. Laurent (1996) builds an efficient and asymptotically normal estimator of $q_2$
 (and more general functionals) for multidimensional distributions using orthogonal projection.
See also references in these papers for more studies under the independence assumption.

In our paper, we study $U$-statistic estimators for $q_{2}$
and $h_{2}$ based on the number of $\ep $-close vector observations (or the number of small inter-point distances)
in a sample from a stationary $m$-dependent sequence with marginal distribution $\mc P$. This extends further the results and
approach in Leonenko and Seleznjev (2010) (see also K\"{a}llberg and Seleznjev,
2012), where the same estimators are studied under independence.
The number of small inter-point distances in an independent sample
exhibits rich asymptotic behaviors, including, e.g., Poisson limits and
asymptotic normality (see Jammalammadaka and Janson, 1986, and references
therein). We show that some of the established limit results for this statistic are
still valid when the sample is from a stationary $m$-dependent sequence.
It should be noted that our normal limit theorems
do not follow from the general theory developed for degenerate variable $U$-statistics under dependence,
see, e.g., Kim et al.,\ 2011, and references therein.

Note that the class of stationary $m$-dependent processes is quite large, see,
e.g., the book of Joe (1997), where there are numerous copula 
constructions for $m$-dependent sequences with given marginal distribution, 
or  Harrelson and Houdre (2003), where the class of stationary
$m$-dependent infinitely divisible sequences is
studied.

First we introduce some notation. Throughout this paper, it is assumed that the
sequence $\{X_{i} \}_{i=1}^{\infty} $ of random $d$-vectors is strictly stationary
and $m$-dependent, i.e., $\{X_{b},X_{b+1}\ldots ,X_{b+s} \}$ and $\{X_{a-r},X_{a-r+1},%
\ldots,X_{a} \}$ are independent sets of vectors when $b-a>m$. Let $\mc{P}$ be
 the (marginal) distribution of $X_t$ with density $p(x), x \in R^d, p(x)\in L_2(R^d)$, and entropy $h_2(\mc{P})$. We write $%
d(x,y):=||x-y||$ for the Euclidean distance in $R^{d}$ and define $%
B_{\ep}(x):=\{y:d(x,y)\leq \ep \}$ to be an $\ep $-ball in $%
R^{d}$ with center at $x$ and radius $\ep $. Denote by $b_{\ep
}(d):=\ep ^{d}b_{1}(d),b_{1}(d)=2\pi ^{d/2}/(d\Gamma (d/2))$, the
volume of the $\ep $-ball. Let $X$ and $Y$ be independent and with
distribution $\mc{P}$ and introduce the $\ep$-ball probability as
\begin{equation*}
p_{X,\ep }(x):=P\{X \in B_{\ep}(x)\}.
\end{equation*}%
Two vectors $x$ and $y$ are said to be $\ep $\textit{-close} if $%
d(x,y)\leq\ep $, for some $\ep >0$. The $\ep $\textit{-coincidence}
probability for independent vectors is written $q_{2,\ep}:=P(d(X,Y)\leq \ep )= \Ex p_{X,\ep }(Y) $.
Then the R\'{e}nyi $\ep $-entropy  $h_{2,\ep }(\mc{P}):=-\log q_{2,\ep }(\mc{P})$
 can be used  as a measure of uncertainty in $\mc{P}$ (see Seleznjev and Thalheim, 2008, Leonenko and  Seleznjev, 2010).
In what follows, let
$\ep =\ep(n) \rightarrow0\mbox{ as }n\rightarrow \infty $. Denote by $|C|$ the
cardinality of the finite set $C$ and let $N_{n }$ be the random
number of $\ep $-close observations in the sample $X_{1},\ldots,X_{n}$,
\begin{equation*}
N_{n}=N_{n,\ep}:=\big |\{d(X_{i},X_{j})\leq \ep ,\;i,j=1,\ldots
,n,\;i<j)\}\big |=\sum_{i<j}I(d(X_{i},X_{j})\leq \ep ) =: \binom{n}{2}Q_n,
\end{equation*}%
where $I(D)$ is the indicator of an event $D$.  Then $Q_n$ is a $U$-statistic of Hoeffding with  varying kernel.
For a short introduction to $U$-statistics techniques, see, e.g., Serfling (2002),
Koroljuk and Borovskich (1994), Lee (1990).
Denote by $\overset{\mathrm{D}}{\rightarrow }$ and $\overset{\mathrm{P}}{%
\rightarrow }$ convergence in distribution and in probability, respectively.
For a sequence of random variables $U_n, n \geq 1$, we write $U_n = \mathrm{O%
}_{\mathrm{P}}(1)$ as $n \to \infty$ if for any $\delta > 0$ and large
enough $n \geq 1$, there exists $C>0$ such that $P(|U_n|>C) \leq \delta$.
Moreover, for a numerical sequence $v_n, n \geq 1$, let $U_n = \Or_{%
\mathrm{P}}(v_n)$ as $n \to \infty$ if $U_n/v_n = \Or_{\mathrm{P}}(1)$
as $n \to \infty$.

The developed technique can also be used for estimation of the
corresponding entropy-type characteristics for discrete distributions (see, e.g., Leonenko and  Seleznjev, 2010) and stationary $m$-dependent sequences. In
this case, the applied estimator is a $U$-statistic with fixed kernel and so the
problem is simplified in the way that some already established general results
yield the limit properties, including consistency and asymptotic normality (see
Appendix). An approach to statistical estimation of the Shannon entropy for %
discrete stationary $m$-dependent sequences can be found in Vatutin and Mikhailov (1995).

The remaining part of the paper is organized as follows.
In Section \ref{sec:Main}, the main results for the number of small inter-point distances $N_{n}$ and the
estimators of $q_{2}$ and $h_{2}$ are presented.
Numerical experiments illustrate the rate of convergence in the obtained asymptotic results.
In Section \ref{sec:App},
we discuss applications of these results to $\ep$-keys in time series databases and
distribution identification problems for dependent samples.
Section \ref{sec:Proofs} contains the proofs of the statements in Section \ref{sec:Main}.
Some asymptotic properties of entropy estimation for the discrete case are given in Appendix.

\section{Main results}

\label{sec:Main}

We formulate the following assumption about finite dimensional distributions
of the stationary sequence $\{X_i\}$.

\begin{enumerate}
\item[ $\mc{A}.$] The marginal density fulfills $p(x) \in L_3(R^d)$.
Moreover, for each $4$-tuple of distinct positive integers $\mathbf{t}%
=(t_1,t_2,t_3, t_4)$, the distribution of the random vector $
(X_{t_1}, X_{t_2},X_{t_3}, X_{t_4})$ has a density $p_{\mathbf{t}}(x_1,x_2,x_3,x_4)$
in $R^{4d}$ that satisfies
\begin{align}  \label{dependence}
g_\mathbf{t}(x_1,x_2):=\left(\int_{R^{2d}} p_{\mathbf{t}}(x_1,x_2,x_3,x_4)^2
dx_3dx_4\right)^{1/2} \in L_1(R^{2d}).
\end{align}
\end{enumerate}

\medskip \noindent \textbf{Remark 1.} (i) The integrability %
\eqref{dependence} ensures that the dependence among $\{X_i \}$ is weak
enough. In fact, provided $p(x) \in L_2(R^d)$, it holds for an independent sequence,
so assumption $\mc{A}$ is a generalization of the
condition $p(x) \in L_3(R^d)$ used for studying the same estimators under
independence (K\"{a}llberg and Seleznjev, 2012).  
\medskip

\noindent 
(ii) If the density $p_\mathbf{t}(x_1,x_2,x_3,x_4)$ is bounded for each distinct $\mathbf{t} = (t_1,t_2,t_3,t_4)$, the following condition is sufficient for $%
\mc{A}$: for each distinct pair $(t_1,t_2)$, let the density $%
p_{t_1,t_2}(x,y)$ of $(X_{t_1},X_{t_2})$ satisfy
\begin{equation*}
p_{t_1,t_2}(x,y) \in L_{1/2}(R^{2d}).
\end{equation*}
\medskip

\noindent Let in the following examples $\{Z_i\}_{i=-\infty}^{\infty}$ be a sequence of independent
identically distributed (i.i.d.) normal $N(0,1)$-random variables.

\medskip \noindent \textbf{Example 1.} (i) Assumption $\mc{A}$ holds
for all vector Gaussian sequences $\{X_i\} $. 
In particular, it is satisfied for the $m$-dependent moving average time
series MA$(m)$ generated by $\{Z_i\}$, i.e.,
\begin{equation}  \label{MA}
X_t=\theta_0 Z_t + \cdots + \theta_mZ_{t-m}, \quad t \ge 1.
\end{equation}

\noindent (ii) An exponential transformation of time series \eqref{MA}
gives a non-linear sequence
\begin{equation*}
X_t = \exp( \theta_0 Z_t + \cdots + \theta_mZ_{t-m}), \quad  t\ge 1.
\end{equation*}
The finite dimensional distributions of $\{X_i\}$ are the multivariate
log-normal distributions (see, e.g., Kotz et al., 2000)
and thus $\mc{A}$ is fulfilled in this case. \medskip

\noindent (iii) Let $X_t =Z_t/Z_{t+1}, t\ge 1$. Then $\{X_i\}$ is a stationary $1$-dependent
 sequence, say, a Cauchy sequence. It can be shown that
for $\mathbf{t}=(1,2,3,4)$,
\begin{align*}
g_\mathbf{t}(x_1,x_2)= C\, \frac{|x_2|}{(1+x_2^2+x_2^2x_1^2)^{5/4}}\in
L_1(R^{2}), \quad C>0,
\end{align*}
that is \eqref{dependence} is valid in this case and
similarly for other $\mathbf{t}$. Since also $p(x) \in L_3(R)$, condition $\mc{A}$ is satisfied.

\subsection{Asymptotic distribution of the number of small inter-point distances}

Let the expectation and variance of the number of small inter-point distances $N_{n}$ be
$\mu_n = \mu_{n,\ep} := \Ex N_{n}$ and $\sigma_n^2 =
\sigma_{n,\ep}^2 := \Var(N_{n})$, respectively.
For $h = 0,1,\ldots$, we
introduce the characteristic $\sigma^2_{1,h,\ep} := \Cov(p_{X,\ep}(X_1),p_{X,\ep}(X_{1+h}))$.
 Let
\begin{equation}\label{zeta1m}
\zeta _{1,m}:=\lim_{n\rightarrow \infty }\frac{1}{n}\Var\left(
\sum_{i=1}^{n}p(X_{i})\right) =\Var(p(X_{1}))+2\sum_{h=1}^{m}\mathrm{%
Cov}(p(X_{1}),p(X_{1+h})).
\end{equation}
\begin{proposition}\label{prop:Nn1}
Suppose that $\mc{A}$ holds.

\begin{itemize}
\item[(i)] Then the expectation and variance of $N_{n}$ fulfill
\begin{align*}
\mu_n & = \binom{n}{2} q_{2,\ep} + \oor(n\ep^{d/2}), \\
\ \sigma_n^2 & = \frac{n^2}{2} q_{2,\ep} + n^3
\left(\sigma^2_{1,0,\ep} + 2\sum_{h=1}^m \sigma^2_{1,h,\ep}
\right) + \oor(n\ep^{d/2}) + \oor(n^2 \ep^{d})
\mbox{
as } n \to \infty.
\end{align*}

\item[(ii)] If $n^2\ep^d \to a, 0<a \leq \infty$, and $\zeta_{1,m}
> 0$ when $\sup_{n \geq 1} \{n\ep^d \} = \infty$, then
\begin{align*}
\mu_n & \sim \frac{1}{2}b_1(d) q_2 n^2\ep^{d}, \\
\sigma_n^2 & \sim \frac{1}{2}b_1(d) q_2 n^2\ep^{d} + b_1(d)^2
\zeta_{1,m} n^3 \ep^{2d} \ninf.
\end{align*}
\end{itemize}
\end{proposition}
\medskip

The asymptotic distribution for $N_{n}$ depends on the rate of decrease of $\ep(n)$.
Some results for $N_n$ under the i.i.d.\ assumption (i.e., $m=0$) are obtained in Jammalamadaka and Janson (1986) (see
also Leonenko and Seleznjev, 2010).
With only additional weak conditions, we show that these results are still valid when $\{X_i\}$ is stationary  and $%
m$-dependent.  Let $\mu = \mu(a) := \frac{1}{2}b_1(d)q_2a$ for $a>0$.

\begin{theorem}\label{th:Nn2} Suppose that $\mc{A}$ holds.

\begin{itemize}
\item[(i)] If $n^2\ep^d \to 0$, then $N_{n} \overset{\mathrm{D}}{\to} 0 %
\ninf$.

\item[(ii)] If $n^2\ep^d \to a$, $0<a<\infty$, then $\mu = \lim_{n\to
\infty} \mu_n$ and
\begin{equation*}
N_{n} \overset{\mathrm{D}}{\to} Po(\mu) \ninf.
\end{equation*}

\item[(iii)] If $n^2\ep^d \to \infty$ and $n\ep^d \to a$, $0 \leq
a \leq \infty$, and $\zeta_{1,m}>0$ when $a=\infty$, then
\begin{equation*}
(N_{n}- \mu_n)/\sigma_n \overset{\mathrm{D}}{\to} N(0,1) \mbox{ as } n \to
\infty.
\end{equation*}
\end{itemize}
\end{theorem}
\medskip

Note that definition \eqref{zeta1m} implies  $\zeta_{1,m} \geq 0$ with equality, e.g., if $\mc P$ is  uniform.

\bigskip 
\noindent \textbf{Remark 2.}
The following inference procedure is discussed in Leonenko and Seleznjev (2010) for i.i.d.\ sequences.
Let $c:=\frac{1}{2}b_1(d)q_2$. By applying
Theorem \ref{th:Nn2}\textit{(ii)} to the minimum inter-point distance {$Y_n
= \min_{1\leq i<j \leq n} ||X_i-X_j||$} and $\ep =
c^{-1/d}t^{1/d}n^{-2/d}$ for a fixed $t>0$, i.e., $\mu_n \to \mu =t$, we
get
\begin{equation*}
P(cn^2Y_n^d > t) =P(Y_n > \ep) = P(N_{n} = 0) \to e^{-\mu} =
e^{-t} \ninf.
\end{equation*}
Hence $Z_n := cn^2Y_n^d$ has asymptotically exponential distribution $Exp(1)$
and an asymptotic confidence interval for the quadratic functional $q_2$ can
be written
\begin{equation*}
I_n = [2c_1/(n^2b_1(d)Y_n^d),2c_2/(n^2 b_1(d) Y^d_n)]
\end{equation*}
for certain positive $c_1,c_2$.

\subsection{Estimation of the entropy-type characteristics}

We consider an estimator for the quadratic functional $q_2$ based on the normalized
statistic $Q_{n}=\binom{n}{2}^{-1} N_{n}$, defined as
\begin{equation*}
\Qt_n = \Qt_{n,\ep} := Q_{n}/\bep.
\end{equation*}
Let $\tilde{H}_n := -\log(\max(\Qt_n,1/n))$ be the corresponding estimator for the entropy $h_2$.
The asymptotic behavior of $\Qt%
_n$ and $\tilde{H}_{n}$ depends on the rate of  decreasing for $\ep(n)$.
In the following theorem for consistency, we give two versions for different asymptotic rates of $\ep(n)$,
with significantly weaker distribution assumptions  in {\it (ii)}.
\begin{theorem}\label{th:P}  \

\begin{itemize}
\item[(i)] If $\mc{A}$ holds and $n^2\ep^d \to \infty$, then
\begin{equation*}
\Qt_n \overset{\mathrm{P}}{\to} q_2 \quad \mbox{and} \quad \tilde{H}_n
\overset{\mathrm{P}}{\to} h_2 \ninf.
\end{equation*}

\item[(ii)] Let $p(x) \in L_3(R^d)$ and assume that $P(X_i \neq X_j)=1$, for
all $i \neq j$. If $n\ep^{d} \to a, 0<a \leq \infty$, then
\begin{equation*}
\Qt_n \overset{\mathrm{P}}{\to} q_2 \quad \mbox{and} \quad \tilde{H}_n
\overset{\mathrm{P}}{\to} h_2 \ninf.
\end{equation*}
\end{itemize}
\end{theorem}
\medskip

Let $\qt_{2,\ep} := q_{2,\ep}/\bep$ and $\tilde{h}%
_{2,\ep} := -\log (\qt_{2,\ep})=  {h}_{2,\ep}+ \log(\bep)$.
Next we show asymptotic normality properties  for the estimators $\Qt_n$ and $\tilde{H}_n$ when $n$ and $\ep$ vary
accordingly. Let $\nu :=~2q_2/b_1(d)$ and recall definition $\eqref{zeta1m}$ of $\zeta_{1,m}$ .

\begin{theorem}\label{th:preLim}
Suppose that $\mc{A}$ holds and $n^2 \ep^d \to \infty$.

\begin{itemize}
\item[(i)] If $n\ep^d \to a, 0<a \leq \infty$, and $\zeta_{1,m}>0$ when
$a=\infty$, then
\begin{equation*}
\sqrt{n} (\Qt_n - \qt_{2,\ep}) \overset{\mathrm{D}}{\to}
N(0, \nu/a+4\zeta_{1,m}) \quad \mbox{and} \quad \sqrt{n} \Qt_{n} (%
\tilde{H}_n - \tilde{h}_{2,\ep}) \overset{\mathrm{D}}{\to} N(0,\nu/a+4\zeta_{1,m}) \ninf.
\end{equation*}

\item[(ii)] If  $n\ep^d \to 0$, then
\begin{equation*}
n\ep^{d/2} (\Qt_n - \qt_{2,\ep}) \overset{\mathrm{D}}{%
\to} N(0, \nu) \quad \mbox{and} \quad n\ep^{d/2} \Qt_{n} (\tilde{H}_n -
\tilde{h}_{2,\ep}) \overset{\mathrm{D}}{\to} N(0, \nu) \mbox{ as } n
\to \infty.
\end{equation*}
\end{itemize}
\end{theorem}
\smallskip

To evaluate the quadratic functional $q_2$ and the entropy $h_2$,
we introduce smoothness conditions for the marginal density $p(x)$.
Denote by $H_{2}^{(\alpha)}(K), 0<\alpha \leq 1, K>0$, a linear space of functions in $R^d$ satisfying a $\alpha$-H\"{o}lder condition in $L_2$-norm with constant $K$, i.e., if $p(x) \in H_2^{(\alpha)}(K)$ and $h \in B_1(d)$, then
\begin{equation}\label{H2}
\left(\int_{R^d} (p(x+h)- p(x))^2 dx \right )^{1/2} \leq K|h|^{\alpha}.
\end{equation}
Note that \eqref{H2} holds, e.g., if for some function $g(x) \in L_2(R^d)$,
\begin{equation}
|p(x+h)-p(x)| \leq g(x) |h|^{\alpha}. \notag
\end{equation}
There are different ways to define the density smoothness, e.g.,
by the conventional or pointwise H\"{o}lder conditions (Leonenko and Seleznjev, 2010, K\"{a}llberg et al., 2012)
or the Fourier characterization (Gin\'{e} and Nickl, 2008).

The rate of convergence in probability can now be described in terms of the
smoothness of $p(x)$.
Let $L(n), n \geq 1$, be a slowly varying function as $n \to \infty$.
\begin{theorem}
\label{th:Cons} Let $\mc{A}$ hold and assume that $p(x) \in
H_2^{(\alpha)}(K) $.

\begin{itemize}
\item[(i)] Then the bias $|\Ex \Qt_n-q_2| \leq \frac{1}{2}K^2
\ep^{2\alpha} + \oor(1/(n\ep^{d/2})) \ninf$.

\item[(ii)] If $0<\alpha\leq d/4$ and $\ep \sim cn^{-2/(4\alpha+d)},
c>0$, then
\begin{equation*}
\Qt_n - q_2 = \Or_{\mathrm{P}}(n^{-\frac{4\alpha}{4\alpha + d}%
}) \quad \mbox{and} \quad \tilde{H}_n - h_2 =\Or_{\mathrm{P}}(n^{-%
\frac{4\alpha}{4\alpha + d}}) \ninf.
\end{equation*}

\item[(iii)] If $\alpha > d/4$ and $\ep \sim L(n)n^{-1/d}$ and $n\ep^d \to a, 0< a \leq\infty$, then
\begin{equation*}
\Qt_n - q_2 = \Or_{\mathrm{P}}(n^{-\frac{1}{2}}) \quad \mbox{and} \quad
\tilde{H}_n - h_2 =\Or_{\mathrm{P}}(n^{-\frac{1}{2}}) \ninf.
\end{equation*}
\end{itemize}
\end{theorem}
\smallskip

\noindent\textbf{Remark 3.} Since the aim of this paper is to provide asymptotic properties for estimation under dependence,
we leave questions regarding efficiency and optimality of the obtained convergence rates for further research.
Nevertheless, it could be mentioned that, in the independent one-dimensional case ($m=0,d=1$), Bickel and Ritov (1988) show that
the rates in Theorem \ref{th:Cons} are optimal in a certain sense (see also Laurent, 1996, Gin\'{e} and Nickl, 2008).
\medskip

In order to make the normal limit results of Theorem \ref{th:preLim}
practical, e.g., to calculate approximate confidence intervals,
the asymptotic variances have to be estimated. In particular, we
need a consistent estimate of the characteristic $\zeta_{1,m}$.
Assuming that $m$ is exactly known might be too strong in our non-parametric setting.
However, note that $\zeta_{1,m} = \zeta_{1,r}$ for $r \geq m$, so 
under the less restrictive assumption that a bound $r$ for $m$ is known, we
can use a consistent estimator of $\zeta_{1,r}$. To construct this estimator, for $h=0,1,\ldots$,
consider the following estimator of $q_{3,h}:=\Ex p(X_1)p(X_{1+h})$,
\begin{equation*}
U_{h,n}=U_{h,n,\ep_0} := M_{h,n}^{-1}b_{\ep_0}(d)^{-2}
\sum_{(i,j,k) \in \mc{E}_{h,n}} I(d(X_{i},X_{j})\leq
\ep_0,d(X_{i+h},X_{k})\leq \ep_0), \quad \ep_0 >0,
\end{equation*}
where $\mc{E}_{h,n}:=\{(i,j,k): 1 \leq i \leq n-(h+1), \; j,k \neq
i,i+h,\; j\ne k\}$ and the number of summands $M_{h,n}:=|\mc{E}_{h,n}| =
(n-(h+1))(n-2)(n-3) $. Let $\ep_0=\ep_0(n) \to 0 \mbox{ as } n
\to \infty$.

\begin{proposition}
\label{std} If $\mc{A}$ holds and $n\ep_0^{3d} \to c, 0<c \leq
\infty$, then
\begin{equation*}
U_{h,n} \overset{\mathrm{P}}{\to} q_{3,h} \ninf.
\end{equation*}
\end{proposition}
\smallskip

\noindent\textbf{Remark 4.} Under the conditions of Proposition \ref{std},
we have
\begin{equation*}
H_{3,n} : = - \frac{1}{2} \log (\max (U_{0,n},1/n)) \overset{\mathrm{P}}{\to}
-\frac{1}{2} \log \left( \int_{R^d} p(x)^3 dx \right) = h_3 \mbox{ as } n
\to \infty,
\end{equation*}
and thus $H_{3,n}$ is a consistent estimator of the {\it cubic} R\'{e}nyi entropy
$h_3$.
\medskip

A consistent plug-in estimator for $\zeta_{1,r}$ can be set up according to
\begin{equation*}
z_{1,r,n}:= U_{0,n}-\Qt_n^2 + 2\sum_{h=1}^r (U_{h,n} - \Qt_n^2),
\end{equation*}
where it is assumed that the sequence $\ep_0= \ep_0(n)$ satisfies $%
n\ep_0^{3d} \to c, 0<c \leq \infty$.

Now we construct asymptotically pivotal quantities by using Theorem \ref{th:preLim}, the smoothness of the marginal density, and
variance estimators. To achieve $\sqrt{n}$%
-rate of convergence, an upper bound $r \geq m$ has to be available.
Let $w^2_{r,n} :=2\Qt_n/(n b_\ep(d)) + 4\max(z_{1,r,n},1/n)$ be the corresponding consistent estimator of
$\nu/a + 4\zeta_{1,m}$ when $n\ep^d \to a, 0<a \leq \infty$.

\begin{theorem}
\label{th:Norm1} Let $\mc{A}$ hold and assume that $p(x) \in
H_2^{(\alpha)}(K), \alpha > d/4$, and $r\geq m$. 
If $\ep \sim L(n)n^{-1/d}$ and $n\ep^d \to a, 0<a \leq \infty$, and $\zeta_{1,m}>0$ when $a=\infty$%
, then
\begin{align*}
\sqrt{n}(\Qt_n- q_2)/w_{r,n} \overset{\mathrm{D}}{\to} N(0,1) \quad %
\mbox{and} \quad \sqrt{n}\Qt_{n}(\tilde{H}_n-h_2)/w_{r,n} \overset{%
\mathrm{D}}{\to} N(0,1) \ninf.
\end{align*}
\end{theorem}
\smallskip

Next we apply Theorem \ref{th:preLim}\textit{(ii)} to weaken the condition $\alpha >d/4$ and therefore get asymptotic normality for less smooth cases. Additionaly, asymptotically pivotal quantities can be built even without a bound $r$ available for $m$.
Note, however, that the obtained rate of convergence is slower than $\sqrt{n}$.
Define a consistent estimator for $\nu$ as $u_n^2:=2\max(\Qt_n,1/n)/b_1(d)$. Let $C_\beta := \beta/(2-\beta), 0 <\beta < 1$.

\begin{theorem}
\label{th:Norm2} Let $\mc{A}$ hold and assume that $p(x) \in
H_2^{(\alpha)}(K)$ and $n^2\ep^d \to \infty$.
\begin{itemize}
\item[(i)] If $\alpha > (d/4)C_\beta$, for some $0<\beta < 1$, and $%
\ep\sim cn^{-(2-\beta)/d}, c>0$, then
\begin{align*}
&n^{\beta/2}c^{d/2}(\Qt_n- q_2)/u_n \overset{\mathrm{D}}{\to} N(0,1)
\quad \mbox{and} \quad n^{\beta/2}c^{d/2}\Qt_n(\tilde{H}_n- h_2)/u_n
\overset{\mathrm{D}}{\to} N(0,1) \ninf.
\end{align*}

\item[(ii)] If $\ep \sim L(n)n^{-2/d}$, then
\begin{align*}
&L(n)(\Qt_n- q_2)/u_n \overset{\mathrm{D}}{\to} N(0,1) \quad \mbox{and}
\quad L(n)\Qt_{n}(\tilde{H}_n-h_2)/u_n \overset{\mathrm{D}}{\to}
N(0,1) \ninf.
\end{align*}
\end{itemize}
\end{theorem}
\smallskip
The practical applicability of the results in this paper relies on an accurate choice of the parameter $\ep$.
One possibility is to use the cross-validation techniques for choosing the optimal bandwidth for density estimation,
see, e.g., Hart and Vieu (1990). However, the problem of finding a suitable $\ep$ is a topic for future research.
\bigskip

\noindent In the following examples, let $\{Z_i\}$ be a sequence of i.i.d. $N(0,1)$-variables.
\medskip

\begin{figure}[hp]
\begin{center}
\includegraphics[width=0.85\textwidth,height=0.35\textheight]{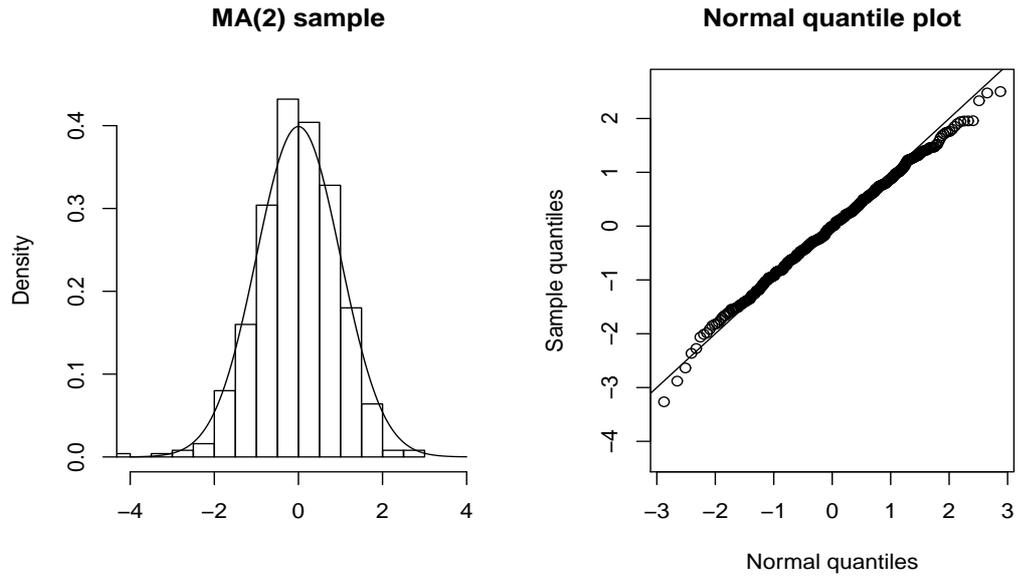}
\end{center}
\caption{A $2$-dependent MA$(2)$ time series;
 sample
size $n=500$, $r=6$, and $\protect\ep = \ep_0=1/10$. Standard normal
approximation for the normalized residuals; $N_{sim} = 500$.}
\label{fig1}
\end{figure}

\noindent \textbf{Example 2.} 
We consider estimation of the quadratic R\'{e}nyi entropy $h_2(\mc{P})$
for the 2-dependent moving average MA$(2)$ process $X_t = \theta_0 Z_{t} +
\theta_1 Z_{t-1}+\theta_2 Z_{t-2}$, where $\theta_2 = \theta_1=\theta_0 = 1/\sqrt{3%
} $. In this case $\mc{P }= N(0,1)$ and $h_2(\mc P) = \log(2\sqrt{\pi})$. We
simulate $N_{sim}=500$ independent and normalized residuals $R^{(i)}_n:=%
\sqrt{n}\Qt_n(\tilde{H}_n-h_2)/w_{r,n}$, $i =1,\ldots,N_{sim}$, with $%
n=500$, $r=6$, and $\ep=\ep_0 =1/10$.
The histogram and normal quantile plot in Figure \ref{fig1} illustrate the
performance of the normal approximation for $R^{(i)}_n$ implied by Theorem~%
\ref{th:Norm1}. The p-value (0.60) for
the Kolmogorov-Smirnov test also supports the hypothesis of standard
normality for the residuals.
\bigskip

\begin{figure}[htbp]
\begin{center}
\includegraphics[width=0.85\textwidth,height=0.35\textheight]{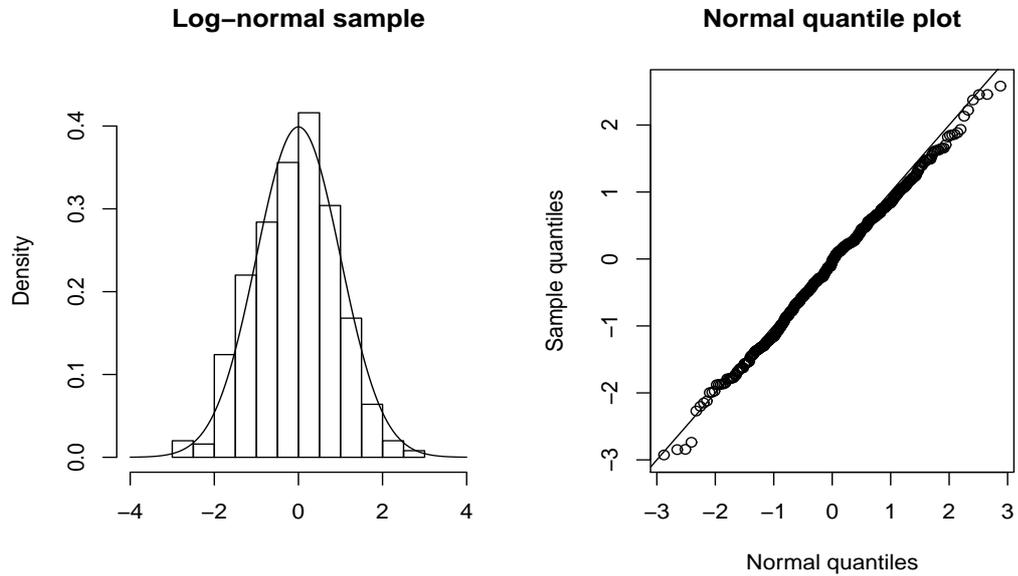}
\end{center}
\caption{A $1$-dependent log-normal sequence; sample size $n=500$, $r=4$,
and $\ep = \ep_0=3/100$. Standard normal approximation for the normalized
residuals; $N_{sim}=500$.}
\label{fig3}
\end{figure}

\noindent \textbf{Example 3.}
Consider the sequence $X_t := \exp(\theta_0 Z_t  +  \theta_1Z_{t-1})$,
where $\theta_0 = \sqrt{3}/2, \theta_{1}=-1/2$, i.e., $\{X_i\}$ is a 1-dependent log-normal sequence,
$X_t \overset{\mathrm{D}}{=} \exp(Z_t)$.
In this case the quadratic entropy $h_2(\mc P) = -\log(e^{1/4}/(2\sqrt{\pi}))$.
Figure \ref{fig3} shows the accuracy of the normal approximation in Theorem~\ref{th:Norm1} for the residuals
$R_{n}^{(i)}:=\sqrt{n}(\Qt_{n}-q_{2})/w_{r,n}$, $i=1,\ldots ,N_{sim}$, where $n=500$, $r=4$, and $\ep = \ep_0 = 3/100$.
The histogram, normal quantile plot, and p-value (0.36) of the Kolmogorov-Smirnov test
imply that the normality hypothesis can not be rejected.
\bigskip

\begin{figure}[hptbp]
\begin{center}
\includegraphics[width=0.85\textwidth,height=0.35\textheight]{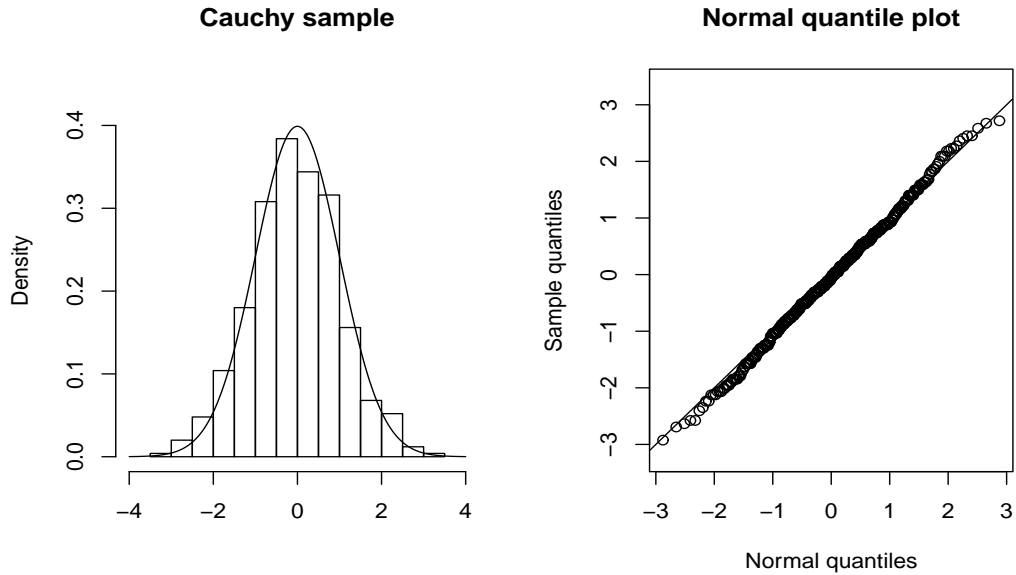}
\end{center}
\caption{ A $1$-dependent Cauchy sequence; sample size $n=500$ and $\protect%
\ep =1/100$. Standard normal approximation for the normalized
residuals; $N_{sim}=500$.}
\label{fig2}
\end{figure}

\noindent \textbf{Example 4.} Estimation of the quadratic functional $q_{2}(%
\mc{P})$ for the 1-dependent Cauchy sequence $X_{t}=Z_{t}/Z_{t+1}$.
Here $\mc{P}$ is the Cauchy distribution with $q_{2}(\mc P)=1/(2\pi) $. We
simulate residuals $R_{n}^{(i)}:=%
n\ep^{d/2}(\Qt_{n}-q_{2})/u_n$, $i=1,\ldots ,N_{sim}$, where $n=500$, $\ep =1/100$, and $N_{sim}=500$.
Figure~\ref{fig2} illustrates the performance
of the normal approximation of $R_{n}^{(i)}$ indicated by Theorem~\ref{th:Norm2}.
The histogram, normal quantile plot, and p-value (0.47) for the
Kolmogorov-Smirnov test allow to accept the hypothesis of standard normality.
\bigskip

\section{Applications}

\label{sec:App}

\noindent \emph{\ $\ep$-Keys in time series databases}

Let a time series {database} $T$ be a matrix with $n$ random records (or
tuples) $r_U(j), j=1,\ldots,n$, and $k$ attributes, $U=\{1,\ldots,k\}$, with
continuous tuple distribution with density $p(x)=p_U(x), x\in R^k$. As
contrast to conventional static databases, the ordering of records in $T$ is
significant, i.e., the timestamp $j$ can be associated with an additional
attribute for $r_U(j)$. For example, time series databases are used for
modelling stock market, environmental (e.g., weather) or web usage data
(see, e.g., Last et al., 2001). Then the database $T$ can be considered as a
sample from a vector time series $\{r_U(t)\}_{t=1}^\infty$ . Assume
additionally that $\{r_U(t)\}$ is a stationary $m$-dependent time series. A
subset $A\subseteq U$, $|A|=d\le k$, is called an $\ep$-\emph{key} if $%
N_{n}(A)=0$, i.e., there are no $\ep$-close in attributes $A$
sub-records $r_A(j),j=1,\ldots,n$. The distribution of $N_{n}(A)$
characterizes the capability of $A$ to distinguish records in $T$ and can be
used to measure the complexity of a database design for further optimization,
e.g., for optimal $\ep$-{key} selection or searching dependencies
between attributes (or \emph{association rules}) (see, e.g., Thalheim,
2000, Seleznjev and Thalheim, 2008, Leonenko and Seleznjev, 2010). Now
Theorem \ref{th:Nn2}\textit{(ii)} gives an approximation of the probability
that $A$ is an $\ep$-key, $P\{N_{n}(A)=0\}\sim e^{-\mu_n}$, where
\begin{equation*}
\mu_n= \frac{n(n-1)}{2} \ep^d b_1(d)\, \qt_{2,\ep} + \mathrm{%
o}(n\ep^{d/2})\sim \frac{1}{2} a b_1(d)\, q_2 \ninf, \;
a>0,
\end{equation*}
i.e., asymptotically optimal $\ep$-key candidates are amongst $A, |A|=d$%
, sets with minimal value of the quadratic functional $q_2$ and the
corresponding estimators of $q_2$ are applicable with various asymptotics
for $\ep$ and $n$ (Remark 2 and Theorems \ref{th:P}, \ref{th:Cons}, \ref%
{th:Norm1}, and \ref{th:Norm2}). \bigskip \\


\noindent \emph{\ Entropy maximizing distributions for  stationary} $m$-\emph{dependent
 sequences}

Note that conditions of consistency for our estimate of the quadratic R\'{e}%
nyi entropy $h_{2}(\mc P)$ (see Theorem \ref{th:P}{\it (ii)}) are rather weak and can
be easily verified for many statistical models. Hence, one can use these
consistent estimators to build goodness of fit tests based on the maximum
entropy principle, see, e.g., Goria et al.\ (2005) and Leonenko and Seleznjev
(2010), where similar approaches were proposed for Shannon and R\'{e}nyi
entropies, respectively. Let us remind some known facts about the maximum
entropy principle, see, e.g., Johnson and Vignat (2007).
Consider the following maximization problem: given a symmetric positive
definite matrix $\Sigma >0$, for all densities $p(x)$ with mean $\mu $
and such that%
\begin{equation}
\int_{\Omega _{0}}p(x)(x-\mu )(x-\mu )^{T}dx=\Sigma ,\quad \Omega
_{0}:=\left\{ x\in R^{d}:(x-\mu )^{T}\Sigma ^{-1}(x-\mu )\leq 4+d\right\} ,
\label{ME0}
\end{equation}%
the quadratic R\'{e}nyi entropy $h_{2}(\mc P)$ is uniquely maximized by the distribution $\mc P^*$ with density
\begin{equation}
p^{\ast }(x):=\left\{
\begin{array}{lr}
A(1-\beta (x-\mu )^{T}\Sigma ^{-1}(x-\mu )), & x\in \Omega _{0}, \\
0, & x\notin \Omega _{0},%
\end{array}%
\right. \;A:=\frac{\Gamma (2+\frac{d}{2})\beta ^{\frac{d}{2}}}{\pi ^{\frac{d%
}{2}}|\Sigma |^{1/2}},\;\beta :=\frac{1}{4+d},  \label{ME1}
\end{equation}%
that is for all other densities support $\Omega _{0}$, mean $\mu$, and
covariance matrix $\Sigma $, see (\ref{ME0}), we have $h_{2}(\mc P)\leq
h_{2}(\mc P^{\ast })$, with equality if and only if $p=p^{\ast }$ almost
everywhere with respect to the Lebesgue measure in $R^{d}$. The distribution
$\mc P^{\ast }$ belongs to the class of multivariate Pearson type II
distributions (or \emph{Student-$r$ distributions}) and its quadratic R\'{e}%
nyi entropy%
\begin{equation}
h_{2}(\mc P^{\ast })=-\log \frac{2\Gamma (2+\frac{d}{2})^{2}\beta ^{\frac{d}{2}%
}}{\Gamma (3+\frac{d}{2})\pi ^{\frac{d}{2}}|\Sigma |^{1/2}}.
\label{ME2}
\end{equation}
For an i.i.d.\ sample, the goodness of fit test based on the maximum quadratic R%
\'{e}nyi  entropy principle was proposed by Leonenko and Seleznjev (2010).
To generalize this test for $m$-dependent data, we need to show that there
exists a stationary $m$-dependent sequence with marginal distribution (\ref%
{ME1}). For one-dimensional processes, one can apply some results from
Joe (1997). Henceforth, we use some definitions and notation from
this book. It is known that for continuous multivariate distributions, the
univariate marginals and the multivariate or dependent structures can be
separated by copula. Let $C(u,v)$ be a bivariate copula with conditional
distribution $C_{\left( 2\right\vert 1)}\left( v\right\vert u)=\partial
C(u,v)/\partial u.$ The inverse conditional distribution is denoted by $%
C_{\left( 2\right\vert 1)}^{-1}\left( s\right\vert u).$\ Let $F$ be a
continuous univariate distribution function and let $\{U _{i}\}$ be a sequence
 i.i.d.\ uniformly distributed $U(0,1)$-random variables. A $1$%
-dependent sequence with marginal distribution $F$ is $Y_{t}=h(U
_{t},U _{t+1})$, where $h(u,v)=F^{-1}[C_{\left(2\right\vert
1)}^{-1}\left( v\right\vert u)]$. The marginal distribution of $Y_{t}$ is $%
F(y)$, see Joe (1997), p.\ 253, and the joint distribution of $(Y_{t},Y_{t+1})
$ is of the form%
\begin{equation*}
P(Y_{t}\leq x,Y_{t+1}\leq y)=\int_{0}^{1}C\left( C_{\left( 2\right\vert
1)}\left( F(x)\right\vert u),F(y)\right) du.
\end{equation*}%
In the same spirit, one can construct a stationary $m$-dependent sequence
with given marginal distribution for any $m\geq 1$, see again Joe (1997),
p.\ 255. In particular, for any $m\geq 1,$ there exists a one-dimensional
stationary $m$-dependent sequence with marginal distribution $F$, which has
density (\ref{ME1}) for $d=1$. For $d\geq 1$, the above copula
construction seems difficult, but one can use again the results from Johnson
and Vignat (2007),  
reformulated for quadratic entropy as follows: if $\{Z_{i}\}$ is a sequence of i.i.d.\ $N(0,1)$-variables, then
\begin{equation*}
X:=\frac{\sqrt{(d+4)}\; \Sigma^{1/2}\cdot Z}{\sqrt{||Z||^{2}+Z_{d+1}^{2}+%
\ldots +Z_{d+4}^{2}}}= f(Z_{1},\ldots ,Z_{d+4}),
\end{equation*}%
has density (\ref{ME1}), where $Z:=(Z_{1},\ldots ,Z_{d})^{T}$. Now a stationary $m$%
-dependent sequence with the marginal density (\ref{ME1}) can be
defined for all $m \in\{1,\ldots, d+3\}$, e.g., $X_{t}:=f(Z_{t},\ldots
,Z_{t+d+3}), t\geq 1, m=d+3$. Consequently, for any fixed $d$ and $m\in
\{1,\ldots, d+3 \}$, there exists a stationary $m$-dependent vector
sequence with marginal maximum quadratic entropy distribution (\ref{ME1}).
Next, from (\ref{ME2}) we get
\begin{equation}\label{ME3}
\frac{e^{h_{2}(\mc P^{\ast })}}{|\Sigma|^{1/2}}=K_d := \frac{\Gamma (3+\frac{d}{2})\pi ^{%
\frac{d}{2}}}{2\Gamma (2+\frac{d}{2})^{2}\beta ^{\frac{d}{2}}}.
\end{equation}

Let $\mc{K}$ be a class of $d$-dimensional density functions $%
p(x),\,x\in R^{d}$, with support $supp\{p\}=\Omega_0$, which satisfy
condition (\ref{ME0}). Note that the density $p^{\ast }(x)$ belongs to this
class. Let $X_{1},\ldots ,X_{n},\,n\geq 2$, be an $m$%
-dependent sample from a member of $\mc{K}$. Consider
$\Sigma _{n}:={1}/{(n-1)}\; \sum_{i=1}^{n}(X_{i}-\bar{X}_{n})(X_{i}-\bar{X}%
_{n})^{T} $, 
the sample covariance matrix, as a consistent estimate of $\Sigma$ as $%
n\rightarrow \infty $, where the sample mean $\bar{X}_{n}:=1/n\;%
\sum_{i=1}^{n}X_{i}$.
The consistency of the sample covariance matrix is a simple consequence of
the asymptotic properties of stationary $m$-dependent processes, see,
e.g., Lee, (1990), Anderson, (1994). Then under the null hypothesis $%
H_{0}:X_{1},\dots X_{n}$ is a sample from\emph{\ }$m$-dependent data
with density $p^{\ast }(x)$, we obtain from the Slutsky theorem and Theorem \ref%
{th:P}{\it (ii)}, that
\begin{equation*}
|\Sigma_{n}|^{-1/2}\exp \left\{ \tilde{H}_{n}\right\} \prob K_{d}\mbox{ as }n\rightarrow \infty ,
\end{equation*}%
where $K_{d}$ is defined in (\ref{ME3}), and $\tilde{H}_{n}$ is the consistent
estimator of the quadratic R\'{e}nyi entropy.

\noindent Under the alternative $H_{1}:X_{1},\dots ,X_{n}$ is a sample
from any other member $p$ of\emph{\ $\mc{K}$}, we find that%
\begin{equation*}
|\Sigma _{n}| ^{-1/2}\exp \left\{ \tilde{H}_{n}\right\} \prob\frac{e^{h_{2}(\mc P)}}{|\Sigma|^{1/2}}<K_{d}\mbox{ as }n\rightarrow
\infty .
\end{equation*}

\noindent In other words, the above mentioned test is consistent against
such alternatives. Note that $K_{1}=\frac{5}{3}\sqrt{5}\simeq
3.\,727,K_{2}=\frac{9}{2}\pi \simeq 14.137,\;K_{3}=\frac{98}{15}\pi \sqrt{%
7}\simeq 54.304$. 
\vspace{0.75cm} 

\section{Proofs}

\label{sec:Proofs} The following lemmas are used in the subsequent proofs.
\noindent

\begin{lemma}[Ch.\ 2, Lee, 1990]
\label{lemma:lee} Let $k \geq 2$. Then the number of k-tuples of integers $1
\leq i_1 < \cdots < i_k \leq n$ that satisfy $i_j-i_{j-1} > h$ for $j =
2,3,\ldots,k$ is $\binom{n-(k-1)h}{k}$.
\end{lemma}

\begin{lemma}[Ch.\ 5, Billingsley, 1995]
\label{Bill95} If, for each k, $Z_n^{(k)} \overset{\mathrm{D}}{\to} Z^{(k)} %
\ninf$, if $Z^{(k)} \overset{\mathrm{D}}{\to} Z$ as $k
\to \infty$, and if
\begin{equation*}
\lim_{k\to \infty} \limsup_{n \to \infty} P(|Z_n^{(k)}-Z_n| \geq \delta) =0,
\end{equation*}
for positive $\delta$, then $Z_n \overset{\mathrm{D}}{\to} Z$.
\end{lemma}


\begin{lemma}
\label{lemma} \

\begin{itemize}
\item[(i)] Assume that $p(x) \in L_3(R^d)$ and let $\pt%
_{X,\ep}(x):=\bep^{-1}p_{X,\ep}(x)$. Then, for $h
=0,1,\ldots$,
\begin{align*}
&\Ex \pt _{X,\ep}(X) =\qt_{2,\ep} \to q_2 \quad \mbox{and} \quad  \Ex \pt%
_{X,\ep}(X_1) \pt_{X,\ep}(X_{1+h}) \to \Ex%
p(X_1)p(X_{1+h})  \mbox{ as }  \ep \to 0,
\end{align*}
and thus $\bep^{-2}\sigma^2_{1,h,\ep} \to \Cov%
(p(X_{1}),p(X_{1+h})) \mbox{ as }  \ep \to 0$.

\item[(ii)] If $\mc{A}$ is satisfied, then
\begin{equation*}
\sup_{i \neq j} P(d(X_i,X_{j}) \leq \ep) = \oor (\ep^{d/2}) %
\mbox{ as } \ep \to 0.
\end{equation*}

\item[(iii)] If $\mc{A}$ is satisfied, then
\begin{equation*}
\sup_{\{i_1,i_2 \} \neq \{j_1,j_2\}} P(d(X_{i_1},X_{j_1}) \leq \ep,
d(X_{i_2},X_{j_2}) \leq \ep) = \oor (\ep^{d}) \mbox{ as }
\ep \to 0.
\end{equation*}
\end{itemize}
\end{lemma}

\medskip \noindent \textit{Proof.} \textit{(i)} 
We use the following result from Lemma 1 in K\"{a}llberg and Seleznjev (2012):
for random vectors $X$ and $Y$ with densities $p_X(x),p_Y(x) \in L_{a+1}(R^d), x \in R^d$, $a \geq 0$,
\begin{equation}\label{eq:pYa}
\Ex \pt_{X,\ep}(Y)^a \to \Ex p(Y)^a \mbox{ as } \ep \to 0.
\end{equation}
Note that \eqref{eq:pYa} immediately implies $\qt_{2,\ep} \to q_2$.
Furthermore, if $Y$ is defined to have density $p_Y(x):=p(x)^2/q_2, x \in R^d$,
we obtain from \eqref{eq:pYa} that
\begin{align}  \label{eq:msCon:Lemma3}
 \Ex\pt_{X,\ep}(X_1)p(X_1) = q_2%
\Ex \pt_{X,\ep}(Y) \to q_2\Ex p(Y) = \Ex p(X_1)^2 \mbox{ as } \ep \to 0.
\end{align}
Now consider the decomposition
\begin{align}  \label{eq:dec:Lemma3}
\Ex \pt_{X,\ep}(X_1)\pt_{X,\ep}(X_{1+h}) =\; &
\Ex p(X_1)p(X_{1+h})+ \Ex  (\pt_{X,\ep}(X_1) -p(X_1) )\pt_{X,\ep}(X_{1+h}) \\
& + \Ex ((\pt_{X,\ep}(X_{1+h}) - p(X_{1+h}))p(X_1).  \notag
\end{align}
By the stationarity of $\{X_i\}$ and H\"{o}lder's inequality, the last two
terms in \eqref{eq:dec:Lemma3} are in absolute value bounded by
\begin{equation*}
\left(\Ex (\pt_{X,\ep}(X_1) - p(X_1) )^2 \right)^{1/2}
\left( \Ex \pt_{X,\ep}(X_{1})^2 \right)^{1/2},
\end{equation*}
which by \eqref{eq:pYa} and \eqref{eq:msCon:Lemma3} tends to zero as $\ep \to 0$.
Hence, $\Ex \pt_{X,\ep}(X_1)\pt_{X,\ep}(X_{1+h}) \to \Ex p(X_1)p(X_{1+h})$ and the assertion follows. \newline
\newline

\noindent \textit{(ii)} First note that the indices $t_3$ and $t_4$ can be
chosen in such a way that $X_{t_3}$ and $X_{t_4}$ are independent and also
independent of $\{X_{t_1},X_{t_2}\}$. For the corresponding density of $%
(X_{t_1},X_{t_2},X_{t_3},X_{t_4})$, we have $p_\mathbf{t}(x_1,x_2,x_3,x_4) =
p_{t_1,t_2}(x_1,x_2)p(x_3)p(x_4)$, where $p_{t_1,t_2}(x_1,x_2)$ is the
density of $(X_{t_1},X_{t_2})$. Consequently, from assumption $\mc{A}$,
\begin{align*}
g_\mathbf{t}(x_1,x_4) & \in L_1(R^{2d}), \\
g_\mathbf{t}(x_1,x_4) & := \left( \int_{R^d}
p_{t_1,t_2}(x_1,x_2)^2p(x_3)^2p(x_4)^2 dx_2dx_3 \right) ^{1/2} = p(x_4)
q_2^{1/2} \left(\int_{R^d} p_{t_1,t_2} (x_1,x_2)^2 dx_2 \right) ^{1/2} .
\end{align*}
Integrating $g_\mathbf{t}(x_1,x_4)$ with respect to $x_4$ gives
\begin{equation*}
g_{t_1,t_2}(x_1) \in L_1(R^d), \quad g_{t_1,t_2}(x_1) := \left( \int_{R^d}
p_{t_1,t_2}(x_1,x_2)^2 dx_2 \right) ^{1/2}.
\end{equation*}
Therefore, by H\"{o}lder's inequality,
\begin{align}  \label{eq1:Lemma3}
P(d(X_{t_1},X_{t_2}) \leq \ep) & = \int_{||x-y||\leq \ep}
p_{t_1,t_2}(x,y) dxdy = \int_{R^d} \left(\int_{||x-y||\leq \ep}
p_{t_1,t_2}(x,y) dy \right) dx \\
& \leq \int_{R^d} \left (\int_{||x-y||\leq \ep} dy \right) ^{1/2} \left (
\int_{||x-y||\leq \ep} p_{t_1,t_2}(x,y)^2 dy \right)^{1/2}dx  \notag \\
& = \bep^{1/2} \int_{R^d} h^{(\ep)}_{t_1,t_2} (x)dx,  \notag
\end{align}
where
\begin{align*}
h^{(\ep)}_{t_1,t_2} (x) :=\left ( \int_{||x-y||\leq \ep}
p_{t_1,t_2}(x,y)^2 dy \right)^{1/2} \leq \left ( \int_{R^d}
p_{t_1,t_2}(x,y)^2 dy \right)^{1/2} = g_{t_1,t_2}(x).
\end{align*}
Since $g_{t_1,t_2}(x) \in L_1(R^d)$ and $h^{(\ep)}_{t_1,t_2} (x) \to 0$
as $\ep \to 0$, the dominated convergence theorem yields $\int_{R^d}
h^{(\ep)}_{t_1,t_2} (x)dx \to 0$ as $\ep \to 0$, and hence from %
\eqref{eq1:Lemma3} we obtain
\begin{equation}  \label{eq3:Lemma3}
P(d(X_{t_1},X_{t_2}) \leq \ep) = \oor(\ep^{d/2}) \mbox{ as }
\ep \to 0,
\end{equation}
for each distinct pair $\{t_1,t_2\}$. Finally, the stationarity and $m$-dependence of the sequence $\{X_i\}$
 imply that $P(d(X_{t_1}, X_{t_2}) \leq \ep)$
attains a finite number of values as $t_1$ and $t_2$ vary.
Thus, the rate of convergence in \eqref{eq3:Lemma3} remains valid if we take the supremum over distinct $%
\{t_1,t_2\}$ as desired, so the statement follows. \newline
\newline

\noindent \textit{(iii)} Since the argument is similar to that of \textit{(ii)},
we show the main steps only. First assume that $i_1=i_2$ and $j_1 \neq
j_2 $. Under $\mc{A}$, it can be shown in a similar way as above that
for each $3$-tuple of distinct positive integers $(t_1,t_2,t_3)$, the
density $p_{t_1,t_2,t_3}(x_1,x_2,x_3)$ of $(X_{t_1},X_{t_2},X_{t_3})$ satisfies
\begin{equation*}
g_{t_1,t_2,t_3}(x_1) \in L_1(R^d), \quad g_{t_1,t_2,t_3}(x_1):= \left(
\int_{R^d}p_{t_1,t_2,t_3}(x_1,x_2,x_3)^2dx_2dx_3 \right)^{1/2}.
\end{equation*}
Now introduce the random vectors $Y_i := (X_{i_1}, X_{i_1})$ and $%
Y_j:=(X_{j_1},X_{j_2})$ in $R^{2d}$ and note that
\begin{equation}  \label{eq4:Lemma3}
P(d(X_{i_1},X_{j_1}) \leq \ep, d(X_{i_1},X_{j_2}) \leq \ep) \leq
P(||Y_i-Y_j||_{2d} \leq 2\ep),
\end{equation}
where $||\cdot||_{2d}$ is the Euclidean norm in $R^{2d}$.
Define the coordinate vectors $y_1:=(x_1,x_1)$ and $y_2:=(x_2,x_3)$.
By H\"{o}lder's inequality,
\begin{align}  \label{eq5:Lemma3}
P(||Y_i-Y_j||_{2d} \leq 2 \ep) & = \int_{||y_1-y_2||_{2d} \leq 2 \ep}
p_{i_1,j_1,j_2}(x_1,x_2,x_3) dx_1 dx_2 dx_3 \\
&\leq b_{2\ep}(2d)^{1/2} \int_{R^{d}} h_{i_1,j_1,j_2}^{(\ep)}(x_1)
dx_1,  \notag
\end{align}
where
\begin{equation}
h_{i_1,j_1,j_2}^{(\ep)}(x_1) := \left(\int_{||y_1-y_2||_{2d} \leq 2
\ep} p_{i_1,j_1,j_2}(x_1,x_2,x_3)^2 dx_2 dx_3 \right)^{1/2} \leq
g_{i_1,j_1.j_2}(x_1).  \notag
\end{equation}
By the dominated convergence theorem,
$\int_{R^d}h_{i_1,j_1,j_2}^{(\ep)}(x_1)dx_1 \to 0$ as $\ep \to 0$
and therefore, since $b_{2\ep}(2d)^{1/2} = C_d \ep^d, C_d >0$, %
we get from \eqref{eq4:Lemma3} and \eqref{eq5:Lemma3} that
\begin{equation}\label{eq:rate2}
P(d(X_{i_1},X_{j_1}) \leq \ep, d(X_{i_2},X_{j_2}) \leq \ep) = \oor%
(\ep^d) \mbox{ as } \ep \to 0,
\end{equation}
where $i_1 = i_2$ and $j_1 \neq j_2$.
By a similar argument, this also valid when $i_1\neq i_2$ and $j_1\neq j_2$.
Finally, from the stationarity and $m$-dependence of $\{X_i\}$,
 the rate of convergence in \eqref{eq:rate2} still holds if we take the supremum over distinct pairs of pairs.
This completes the proof. \hfill $\Box$ \newline
\newline

\noindent \textit{Proof of Proposition \ref{prop:Nn1}.} \textit{(i)} Define
the index set $\mc{I }= \mc{I}({n,m}) := \{(i,j): 1\leq i < j \leq
n, j-i >m \}$ and the reduced form of $N_{n}$:
\begin{equation*}
N^*_{n} := \sum_{(i,j) \in \mc{I}} I(d(X_i,X_j)\leq \ep).
\end{equation*}
 Since $\Ex I(d(X_i,X_j)\leq \ep) = q_{2,\ep}$ when
$(i,j) \in \mc{I}$, Lemmas \ref{lemma:lee} and \ref{lemma} yield the
claim for the expectation:
\begin{align}  \label{eq:ENn}
\Ex N_{n} & = \Ex N^*_{n} + \sum_{h=1}^m (n-h)
P(d(X_1,X_{1+h})\leq \ep) = \binom{n-m}{2}q_{2,\ep} + \oor%
(n\ep^{d/2}) \\
& = \binom{n}{2}q_{2,\ep} + \oor(n\ep^{d/2}) \mbox{ as } n
\to \infty.  \notag
\end{align}
For the variance of $N_{n}$, we first study the variance of $N^*_{n}$. To
this end, we need to consider $|\mc{I}|^2=\binom{n-m}{2}^2$ terms of
the form
\begin{equation}  \label{cov}
\Cov(I(d(X_{s_1},X_{s_2})\leq \ep),I(d(X_{t_1},X_{t_2})\leq \ep)),
\end{equation}
where $(s_1,s_2), (t_1,t_2) \in \mc{I}$. To count the various types of
such terms, we use results from Ch.\ 2.4.1,  Theorem 1,  in Lee~(1990). It should be noted
that these results are stated for $U$-statistics based on random variables,
rather than random vectors, but the argument is merely combinatorial and
thus also valid here.

\begin{enumerate}
\item[1)] When $|s_i - t_j| > m, i,j=1,2$, the random variables in $%
\eqref{cov}$ are independent and the covariance is zero.

\item[2)] Only two random variables are involved, i.e., $s_1 = t_1$ and $%
s_2=t_2$, so the covariance is $\Var(I(d(X_{s_1},X_{s_2})<
\ep)) = q_{2,\ep} - q^2_{2,\ep} = q_{2,\ep} + \Or%
(\ep^{2d})$ as $\ep \to 0$. The number of such terms is $|\mc{%
I}| = \binom{n-m}{2}$.

\item[3)] Exactly one of the four possible differences $|s_i-t_j|$ is zero
and the rest are greater than $m$. By conditioning, the covariance is $%
\Cov(I(d(X_{s_1},X_{s_2})\leq \ep),I(d(X_{s_1},X_{t_2})\leq \ep))
= \sigma^2_{1,0,\ep}$. There are $6\binom{n-2m}{3}$ terms of this type.

\item[4)] For $h=1,\ldots,m$, $0 < |s_i-t_j| = h \leq m$ for one of the
differences $|s_i-t_j|$ and the others are greater than $m$.
Then the covariance is $\Cov%
(I(d(X_{s_1},X_{s_2})\leq \ep),I(d(X_{s_1+h},X_{t_2})\leq \ep)) =
\sigma^2_{1,h,\ep}$. The number of terms of this type is $12\binom{%
n-2m-h}{3}$.

\item[5)] The number of the remaining terms is $\Or(n^2)$, so Lemma %
\ref{lemma}\textit{(iii)} indicates that their sum is $\oor%
(n^2\ep^{d}) \ninf$.
\end{enumerate}
From the above and Lemma \ref{lemma}, we obtain
\begin{align}  \label{VarN*}
\Var(N^*_{n}) = & \sum_{(s_1,s_2) \in \mc{I}} \sum_{(t_1,t_2)
\in \mc{I}} \Cov(I(d(X_{s_1},X_{s_2})\leq
\ep),I(d(X_{t_1},X_{t_2})\leq \ep)) \\
=& \binom{n-m}{2}q_{2,\ep} + 6\binom{n-2m}{3}\sigma^2_{1,0,%
\ep} + \sum_{h=1}^m 12 \binom{n-2m-h}{3}\sigma^2_{1,h,\ep}+
\oor(n^2\ep^{d})  \notag \\
=& \frac{n^2}{2} q_{2,\ep} + n^3 \left(\sigma^2_{1,0,\ep} +
2\sum_{h=1}^m \sigma^2_{1,h,\ep} \right)+ \oor(n^2\ep^{d}) %
\ninf.  \notag
\end{align}
Next it is verified that $\Var(N^*_{n})$ and $\Var(N_{n})$
are close. Note that
\begin{equation}  \label{VarN}
\Var(N_{n}) = (n-1)C_{1,n}+\ldots+(n-m)C_{m,n}+\Cov%
(N_{n},N^*_{n}),
\end{equation}
where $C_{h,n} := \Cov(N_{n}, I(d(X_1,X_{1+h})\leq \ep)), h=
1,\ldots,m$. We have
\begin{align*}
C_{h,n}& = \sum_{1\leq i<j \leq n} \Cov(I(d(X_i,X_j)\leq
\ep),I(d(X_1,X_{1+h})\leq \ep) \\
& = \Var(I(d(X_1,X_{1+h})\leq \ep)) + \sum_{\substack{ 1\leq i<j
\leq n  \\ \{i,j\} \neq \{1,1+h\}}} \Cov(I(d(X_i,X_j)\leq
\ep),I(d(X_1,X_{1+h})\leq \ep)),
\end{align*}
where the number of pairs $(i,j)$ in the last sum for which $I(d(X_i,X_j)\leq
\ep)$ and $I(d(X_1,X_{1+h})\leq \ep)$ are \textit{not} independent is
$\Or(n)$. From Lemma \ref{lemma}, we get $C_{h,n} = \oor%
(\ep^{d/2}) + \oor(n \ep^{d}) \ninf$,
and consequently \eqref{VarN} leads to
\begin{equation}  \label{VarN1}
\Var(N_{n}) = \Cov(N_{n},N^*_{n}) + \oor%
(n\ep^{d/2}) + \oor(n^2\ep^{d}) \ninf.
\end{equation}
Furthermore, observe that
\begin{equation*}
\Var(N_{n}) = \Var(N^*_{n}) + \Var(N_{n} - N^*_{n})
+ 2\Cov(N^*_{n},N_{n}-N^*_{n})
\end{equation*}
and hence \eqref{VarN1} gives
\begin{equation}  \label{VarN2}
\Var(N_{n}) = \Var(N^*_{n}) - \Var(N_{n}- N^*_{n}) +
\oor(n\ep^{d/2}) + \oor(n^2\ep^{d}) \mbox{ as } n \to
\infty.
\end{equation}
The next step is to show that
\begin{equation}\label{V(N*-N)}
\Var(N_{n}-N^*_{n}) = \oor(n\ep^{d/2})
+ \oor(n^2\ep^{d}) \ninf.
\end{equation}
To this end, first note that
Lemma \ref{lemma}\textit{(ii)} implies
\begin{equation}\label{eq:exConv}
\Ex(N_{n} - N^*_{n}) = \sum_{\substack{ 1\leq i<j \leq n  \\ j-i \leq
m}} \Ex I(d(X_i,X_j)\leq \ep) = \oor(n\ep^{d/2}) \mbox{
as } n \to \infty,
\end{equation}
since the number of terms in this sum is $\Or(n)$, so we
only need to prove that
\begin{equation}  \label{eq:msConv}
\Ex(N_{n}-N^*_{n})^2 = \oor(n\ep^{d/2}) + \oor%
(n^2\ep^{d}) \ninf.
\end{equation}
We have
\begin{equation*}
\Ex(N_{n}-N^*_{n})^2 = \sum_{\substack{ 1\leq i<j \leq n  \\ j-i \leq
m}}\sum_{\substack{ 1\leq i^{\prime}<j^{\prime}\leq n  \\ %
j^{\prime}-i^{\prime}\leq m}} \Ex I(d(X_i,X_j)\leq
\ep)I(d(X_{i^{\prime}},X_{j^{\prime}})\leq \ep).
\end{equation*}
There are $\Or(n^2)$ terms in this double sum. Moreover, $\Or%
(n)$ of these have the form $P(d(X_i,X_{j})\leq \ep)$ and the remaining are
of the type $P(d(X_{i_1},X_{j_1}) \leq \ep,d(X_{i_2},X_{j_2})\leq \ep),$
where $\{i_1,j_1\} \neq \{i_2,j_2\}$,
so \eqref{eq:msConv} follows from Lemma 3.
Therefore \eqref{V(N*-N)} holds true, and combining  with \eqref{VarN2} and \eqref{VarN*}, we get
\begin{align}  \label{VarNn}
\Var(N_{n})= \frac{n^2}{2} q_{2,\ep} + n^3
\left(\sigma^2_{1,0,\ep} + 2\sum_{h=1}^m \sigma^2_{1,h,\ep}
\right) + \oor(n\ep^{d/2})+\oor (n^2\ep^{d})
\mbox{ as
} n \to \infty.
\end{align}
The assertion is proved.\newline
\newline
\textit{(ii)} If $\zeta_{1,m} > 0$, the claim follows directly from Lemma %
\ref{lemma} and \eqref{VarNn}. When $\zeta_{1,m}=0$, we have $\sup _{n \geq 1}\{n\ep^d \} < \infty$ by assumption.
Hence, from  \eqref{VarNn} and the condition $n^2\ep^d \to a, 0<a \leq \infty,$
\begin{equation*}
\frac{\Var (N_{n})}{\frac{1}{2}b_1(d)q_2n^2\ep^d} = 1 + \oor%
(n\ep^d) + \oor(1/(n\ep^{d/2})) + \oor(1) \to 1 \mbox{
as } n \to \infty.
\end{equation*}
This completes the proof. \hfill $\Box$ \newline
\newline

\noindent \textit{Proof of Theorem \ref{th:Nn2}.} \textit{(i)} From the
assumption $n^2\ep^d \to 0$ and Proposition \ref{prop:Nn1}\textit{(i)},
we get $\mu_n$, $\sigma_n^2 \to 0$. Consequently,  $N_{n} \overset{%
\mathrm{P}}{\to} 0$ and the assertion follows.\newline
\newline
\textit{(ii)} By the condition $n^2\ep^d \to a, 0<a < \infty$, we have $%
n\ep^d \to 0$ and thus Proposition \ref{prop:Nn1}\textit{(i)} and
Lemma~\ref{lemma}\textit{(i)} yield $\mu_n \to \mu$. Moreover,
since \eqref{eq:msConv} and $n^2\ep^d \to a, 0<a < \infty$, imply $\Ex(N_{n}^* - N_{n})^2 \to 0$,
it is enough to verify the Poisson convergence
for $N^*_{n}$. For this we apply the Stein-Chen Poisson approximation method
(Barbour et al., 1992). To measure the deviation between two probability
distributions $\mc{P}_1$ and $\mc{P}_2$ of integer-valued random
variables, we use the \textit{total variation distance}
\begin{equation*}
d_{TV}(\mc{P}_1, \mc{P}_2) := \sup\{|\mc{P}_1(C)-\mc{P}%
_2(C)|:C \subset \{0,1,\ldots\}\} = \frac{1}{2} \sum_{i \geq 0} |\mc{P}%
_1(i) - \mc{P}_2(i)|.
\end{equation*}
Recall the definition $\mc{I }:= \{(i,j): 1\leq i < j \leq n, j-i >m \}$ and let
$\mc{I}_{(i,j)} := \mc{I }\setminus \{(i,j)\}$.
Further, using the notation $I_\ep(i,j):=I(d(X_i,X_j)\leq \ep)$, we define
\begin{equation*}
\mc{I}_{(i,j)}^{0} := \{(k,l) \in \mc{I}%
_{(i,j)}: I_\ep(k,l) \mbox{ and } I_\ep(i,j) \mbox{ are dependent}\}.
\end{equation*}
Since $\Ex I_\ep(i,j) = q_{2,\ep}$ for all $(i,j) \in \mc{I}
$, Corollary 2.C.5 in Barbour et al.\ (1992) yields
\begin{align}  \label{eq:tdv}
d_{TV}(\mc{L}(N^*_{n}), & Po(\mu^*_{n})) \\
&\leq \frac{1}{\mu^*_{n}} \Bigg(\sum_{(i,j) \in \mc{I}}
q^2_{2,\ep} + \sum_{(i,j) \in \mc{I}}\sum_{(k,l) \in \mc{I}%
_{(i,j)}^{0}} q^2_{2,\ep} + \sum_{(i,j) \in \mc{I}}\sum_{(k,l) \in
\mc{I}_{(i,j)}^{0}} \Ex I_\ep(i,j)I_\ep(k,l) \Bigg),  \notag
\end{align}
where $\mc{L}(N^*_{n})$ denotes the distribution of $N^*_{n}$ and $%
\mu^*_{n}:= \Ex N^*_{n}$. The $m$-dependence of $\{X_i\}$ gives
that the number of terms in $\mc{I}^0_{(i,j)}$, for each $(i,j) \in
\mc{I}$, has the bound
\begin{equation*}
|\mc{I}^0_{(i,j)}| \leq n(4m+2),
\end{equation*}
so for the first two sums in \eqref{eq:tdv}, Lemma \ref{lemma}
and the condition $n^2\ep^d \to a, 0<a<\infty$, lead to
\begin{equation}  \label{eq:poiss:bound1}
\sum_{(i,j) \in \mc{I}} q^2_{2,\ep} \leq \sum_{(i,j) \in \mc{I%
}}\sum_{(k,l) \in \mc{I}_{(i,j)}^{0}} q^2_{2,\ep} \leq \binom{n-m}{%
2}n(4m+2) q^2_{2,\ep} \sim (2m + 1)b_1(d)^2q_2^2 a^2n^{-1} \mbox{ as
} n \to \infty.
\end{equation}
For the last sum in \eqref{eq:tdv}, we see that the various types of terms
can be counted in a similar way as the covariances in \eqref{VarN*}.
With the notation $q_{3,h,\ep} := \Ex p_{X,\ep}(X_{1})p_{X,\ep}(X_{1+h})$,
it follows from $n^2\ep^d \to a,
0<a<\infty$, and Lemma \ref{lemma}\textit{(i)} that
\begin{align}  \label{eq:poiss:bound2}
\sum_{(i,j) \in \mc{I}} &\sum_{(k,l) \in \mc{I}_{(i,j)}^{0}}
\Ex I_\ep(i,j)I_\ep(k,l) = 6\binom{n-2m}{3} q_{3,0,\ep}+%
\sum_{h=1}^m 12 \binom{n-2m-h}{3} q_{3,h,\ep} + \oor%
(n^2\ep^d) \\
& = b_1(d)^2a^2 \left(\Ex p(X_1)^2 +2\sum_{h=1}^m \Ex
p(X_1)p(X_{1+h}) \right)n^{-1} +\oor(n^{-1}) + \oor%
(1) \to 0 \ninf.  \notag
\end{align}
Note also that
\begin{equation}  \label{eq:mu^*}
\mu^*_{n} = \binom{n-m}{2} q_{2,\ep} \to \mu \ninf.
\end{equation}
By combining \eqref{eq:tdv}, \eqref{eq:poiss:bound1}, \eqref{eq:poiss:bound2}, and \eqref{eq:mu^*},
\begin{equation}  \label{eq:poiss:conv}
d_{TV}(\mc{L}(N^*_{n}), Po(\mu^*_{n})) \to 0 \mbox{
as } n \to \infty.
\end{equation}
Finally, if $Z_n \sim Po(\mu^*_{n})$ and $Z \sim Po(\mu)$, then \eqref{eq:mu^*} yields $Z_n \overset%
{\mathrm{D}}{\to} Z$, so from \eqref{eq:poiss:conv} we obtain
\begin{equation*}
P(N^*_{n}=k) - P(Z=k) = P(N^*_{n}=k)-P(Z_n=k) + P(Z_n=k)-P(Z=k) \to 0 \mbox{
as } n \to \infty,
\end{equation*}
for each $k=0,1,\ldots$, and thus $N^*_{n} \overset{\mathrm{D}}{\to} Z$.
The statement follows. \newline
\newline
\textit{(iii)} The idea of the proof is to apply Lemma \ref{Bill95} to the
reduced form $N^*_{n}$ of $N_{n}$. In fact, by \eqref{eq:exConv}, \eqref{eq:msConv}, and Proposition \ref{prop:Nn1}{\it (ii)},
for the last two terms in the decomposition
\begin{equation*}
\frac{N_{n} - \Ex N_{n}}{\sigma_n} = \frac{N^*_{n}-\Ex N^*_{n}%
}{\sigma_n} + \frac{\Ex (N^*_{n} - N_{n})}{\sigma_n}+ \frac{%
N_{n}-N^*_{n}}{\sigma_n} ,
\end{equation*}
we have ${\Ex (N^*_{n} - N_{n})}/{\sigma_n}\to 0 $ and
$({N_{n}-N^*_{n}})/{\sigma_n}\to 0$ (in quadratic mean).
Consequently,
it is enough to show
\begin{equation}  \label{eq:Zn}
Z_n  := \frac{N^*_{n}- \Ex N^*_{n}}{\sigma_n} \overset%
{\mathrm{D}}{\to} N(0,1) \ninf.
\end{equation}
The weak convergence \eqref{eq:Zn} follows from Lemma \ref{Bill95} if
there exist successive approximations $\{Z_n^{(k)}\}_{k=1}^\infty$ of $Z_n$
such that

\begin{enumerate}
\item[c1)] $Z^{(k)}_n \overset{\mathrm{D}}{\to} Z^{(k)} \mbox{ as } n \to
\infty$ and $Z^{(k)} \overset{\mathrm{D}}{\to} N(0,1) \mbox{ as } k \to
\infty$.

\item[c2)] For every $\delta > 0$,
\begin{equation*}
\lim_{k \to \infty} \limsup_{n \to \infty} P(|Z_n-Z_n^{(k)}| \geq \delta) =0.
\end{equation*}
\end{enumerate}
In order to construct such $Z_n^{(k)}$, for each $k= 1,2, \ldots$, define an
integer $s = s(k,n) := [n/(k+m)]$ and consider the set of $k$-subsets of $%
\{1,\ldots,n\}$
\begin{equation*}
S^{(k)}_i := \{(i-1)(k+m)+1,\ldots,(i-1)(k+m)+k\}, \quad i = 1, \ldots, s.
\end{equation*}
Now, based on the subset $\mc{I}^{(k)}:= \{(i,j): i \in S_l^{(k)}, j
\in S_t^{(k)}, 1\leq l<t \leq s\}$ of $\mc{I}$, let
\begin{equation*}
Z^{(k)}_n := \frac{U^{(k)}_n - \Ex U^{(k)}_n}{\sigma_n}, \quad
U^{(k)}_{n} = U^{(k)}_{n,\ep} := \sum_{(i,j)\in\mc{I}^{(k)}}
I(d(X_i,X_j)\leq \ep).
\end{equation*}

First we prove c2).
Denote by $M_2$, $M_{3}$, and $\{M_{4,h}\}_{h=1}^m$, the numbers of terms of
types 2, 3, and 4, respectively, defined in Proposition \ref{prop:Nn1}%
\textit{(i)}. Furthermore, let $M^{(k)}_2$, $M^{(k)}_{3}$, and $%
\{M^{(k)}_{4,h}\}_{h=1}^m$ be the numbers of these terms that also appear in
\begin{equation*}
\Var(U_n^{(k)}) = \sum_{(s_1,s_2)\in \mc{I}^{(k)}}
\sum_{(t_1,t_2)\in \mc{I}^{(k)}} \Cov(I(d(X_{s_1},X_{s_2})\leq
\ep),I(d(X_{t_1},X_{t_2})\leq \ep)).
\end{equation*}
Since $\mc{I}^{(k)} \subseteq \mc{I}$, we observe that
\begin{align}  \label{VarNU}
\Var(N^*_{n} - U_n^{(k)}) \leq &\sum_{(s_1,s_2) \in \mc{I}%
}\sum_{(t_1,t_2) \in \mc{I}} |\Cov(I(d(X_{s_1},X_{s_2})\leq
\ep),I(d(X_{t_1},X_{t_2})\leq \ep))| \\
& -\sum_{(s_1,s_2) \in \mc{I}^{(k)}}\sum_{(t_1,t_2) \in \mc{I}%
^{(k)}} |\Cov(I(d(X_{s_1},X_{s_2})\leq \ep),I(d(X_{t_1},X_{t_2})\leq
\ep))|  \notag \\
= & \; (M_2- M_2^{(k)})q_{2,\ep} +
(M_3-M^{(k)}_3)\sigma^2_{1,0,\ep} + \sum_{h=1}^m
(M_{4,h}-M^{(k)}_{4,h})|\sigma^2_{1,h,\ep}|  \notag \\
& + \oor(n^2\ep^{d}) \ninf.  \notag
\end{align}
First, for $h = 1,\ldots,m$, we obtain a lower bound for $M^{(k)}_{4,h}$.
Assume without loss of generality that $k \geq 2m+1$.
There are $\binom{s}{2}(k-2h)^2$ elements $(s_1,s_2)$ in $\mc{I}^{(k)} $ that
also satisfy $(s_1 + d_1, s_2 + d_2) \in \mc{I}^{(k)}$, where $d_1,d_2 =
\pm h$. For every element of this type, an index $t_1$ with $|s_1 - t_1| = h$
or $|s_2-t_1|=h$ can be chosen in 4 different ways. Moreover, for each such
alternative we can choose $t_2$ in at least $(s-2)k$ different ways.
We get that $M_{4,h}^{(k)}$ is bounded from below by
\begin{equation}\label{M4h}
\bar{M}_{4,h}^{(k)}:=\binom{s}{2}(k- 2h)^24(s-2)k \sim 2 \frac{(k-2h)^2 k}{(k+m)^3} n^3 \ninf.
\end{equation}
Further, by Proposition \ref{prop:Nn1}\textit{(ii)}, Lemma \ref{lemma}{\it (i)}, and the limits $%
M_{4,h} \sim 2 n^3$ and \eqref{M4h}, we get
\begin{align}  \label{M2h}
\frac{(M_{4,h}-\bar{M}^{(k)}_{4,h})|\sigma^2_{1,h,\ep}|}{\sigma^2_n} &
= \frac{n^3\ep^{2d}}{\sigma^2_n}\frac{(M_{4,h}-\bar{M}%
^{(k)}_{4,h})|\sigma^2_{1,h,\ep}|}{n^3\ep^{2d}} \\
&\to C_{4,h} \left(1-\frac{(k-2h)^2 k}{(k+m)^3} \right)|\Cov%
(p(X_1),p(X_{1+h}))| \ninf,  \notag
\end{align}
for some $0 \leq C_{4,h} < \infty$.
By a similar argument, we have lower bounds $\bar{M}_2^{(k)}$ and $\bar{M}%
_3^{(k)}$ for $M_2^{(k)}$ and $M_3^{(k)}$, respectively, such that, for some $0
\leq C_2 < 0$ and $0 \leq C_3 < \infty$,
\begin{align}  \label{M13}
\frac{(M_2- \bar{M}_2^{(k)})q_{2,\ep}}{\sigma^2_n} = \frac{n^2\ep^d%
}{\sigma^2_n}\frac{(M_2- \bar{M}_2^{(k)})q_{2,\ep}}{n^2\ep^d} &\to
C_2\left(1 - \frac{k^2}{(k+m)^2}\right) q_2, \\
\frac{(M_3- \bar{M}_3^{(k)})\sigma^2_{1,0,\ep}}{\sigma_n^2} = \frac{%
n^3\ep^{2d}}{\sigma_n^2}\frac{(M_3- \bar{M}_3^{(k)})\sigma^2_{1,0,%
\ep}}{n^3\ep^{2d}} & \to C_3 \left(1 - \frac{k^3}{(k+m)^3} \right)
\Var(p(X_1)) \ninf.  \notag
\end{align}
Now, from \eqref{VarNU}, \eqref{M2h}, \eqref{M13}, and Proposition \ref{prop:Nn1}{\it (ii)},
\begin{equation*}
\lim_{k\to \infty} \lim_{n\to \infty} \Var (Z_n-Z_n^{(k)}) =
\lim_{k\to \infty} \lim_{n\to \infty} \Var \left(\frac{%
N^*_{n}-U_n^{(k)}}{\sigma_n} \right) = 0
\end{equation*}
and hence, since $Z_n-Z_n^{(k)}$ has zero mean, c2) is implied by
Chebyshev's inequality.

Next we prove that c1). Let
\begin{equation*}
V(k):= \left \{
\begin{array}{lr}
\dfrac{\dfrac{1}{2}\dfrac{k^2}{(k+m)^2} q_2 + \dfrac{k^3}{(k+m)^3}%
b_1(d)\zeta^{(k)}_{1,m} a}{\dfrac{1}{2} q_2+ b_1(d)\zeta_{1,m}a}, & \quad 0 \leq a
< \infty, \\ [2em]
\dfrac{k^3 }{(k+m)^3 } \dfrac{\zeta^{(k)}_{1,m}}{\zeta_{1,m}}, & \quad a = \infty,%
\end{array}
\right.
\end{equation*}
where
$$
\zeta_{1,m}^{(k)} := \Var(p(X_1)) + 2 \sum_{h=1}^m (1- h/k) \Cov(p(X_1),p(X_{1+h})).
$$
Note that, since $\zeta_{1,m}^{(k)} \to \zeta_{1,m}$ as $k \to \infty$, we have $V(k) \to 1$ as $k \to \infty$, so c1) follows if it can be verified that
\begin{equation}\label{N(0,V(k))}
Z_n^{(k)} \to N(0,V(k)) \ninf.
\end{equation}
To prove this, we apply the corresponding result of Jammalamadaka and Janson (1986) for independent samples.
In fact, if we introduce the pooled random vectors in $R^{kd}$
\begin{equation*}
Y_i := (X_{(i-1)(k+m)+1},\ldots, X_{(i-1)(k+m)+k}), \quad i= 1, \ldots, s,
\end{equation*}
then the $m$-dependence of $\{X_i \}$ implies that $\{Y_i\}$ is an \textit{%
independent} sequence in $R^{kd}$. Thus $U_n^{(k)}$ can be represented as a
 $U$-statistic with respect to the independent sample $Y_1,\ldots, Y_s$,
\begin{equation*}
W_s^{(k)}:=U_n^{(k)} = \sum_{1 \leq i<j \leq s} f^{(k)}_s(Y_i,Y_j), \quad
f_s^{(k)}(Y_i,Y_j) := \sum_{l \in S_i^{(k)}} \sum_{t \in S_j^{(k)}}
I(d(X_l,X_t)\leq \ep).
\end{equation*}
Furthermore, let $g_s^{(k)}(Y_1) := \Ex (f_s(Y_1,Y_2)|Y_1) =
k\sum_{i=1}^k p_{X,\ep}(X_i)$ and define
\begin{equation} \label{eta1}
\eta^2_{s,k}:= \frac{1}{2}s^2 \Var(f_s^{(k)}(Y_1,Y_2)) + s^3\mathrm{%
Var}(g_s^{(k)}(Y_1)).
\end{equation}
Since $n^2\ep^d \to \infty$ yields $s^2 \ep^d \to \infty$, we obtain from Lemma \ref{lemma} that
\begin{equation}  \label{eta}
\eta^2_{s,k} \geq \frac{1}{2}s^2\Var(f_s^{(k)}(Y_1,Y_2))  \sim \frac{1}{2}b_1(d)q_2 (sk)^2 \ep^d \to \infty \ninf,
\end{equation}
and therefore
\begin{equation}  \label{c1}
\sup_{y_i,y_j} |f^{(k)}_s(y_i,y_j)| = k^2 = \oor(\eta_{s,k})
\ninf.
\end{equation}
Moreover, Jammalamadaka and Janson (1986) show that $\sup_x
p_{X,\ep} (x) = \oor(\ep^{d/2})$ as $\ep \to 0$, and hence the stationarity of $\{X_i\}$ and \eqref{eta} give
\begin{equation}  \label{c2}
\sup_{y} \Ex |f^{(k)}_s(y,Y_1)|  \leq k \sum_{i=1}^k \sup_x \Ex
I(d(x,X_i)\leq \ep) = k^2 \sup_x p_{X,\ep} (x) = \oor(\ep^{d/2}) =
\oor(\eta_{s,k}/s) \ninf.  \notag
\end{equation}
Since $n \to \infty$ implies $s \to \infty$, we get from \eqref{c1} and \eqref{c2} that the conditions of Theorem 2.1 in
Jammalamadaka and Janson (1986) are satisfied. Consequently
\begin{equation}\label{Znk}
\frac{\sigma_n}{\eta_{s,k}} Z_n^{(k)}=\frac{W_s^{(k)}- \Ex W_s^{(k)}}{\eta_{s,k}}   \overset{\mathrm{D}} {\to}
N(0,1) \ninf.
\end{equation}
Further, using Proposition \ref{prop:Nn1}\textit{(ii)} and  definition \eqref{eta1},
it is straightforward to show
$$\eta_{s,k}^2/\sigma_n^2 \to V(k) \ninf,$$
so the desired limit \eqref{N(0,V(k))} of $Z_n^{(k)}$ follows from \eqref{Znk} and the Slutsky theorem.
This completes the proof. \hfill $\Box$ \newline
\newline

\noindent \textit{Proof of Theorem \ref{th:P}.} \textit{(i)} From Lemma \ref%
{lemma}, Proposition \ref{prop:Nn1}\textit{(i)}, and the condition $n^2
\ep^d \to \infty$,
\begin{align}  \label{eq:VarQtn}
\Ex\Qt_n & = \qt_{2,\ep} + \oor\left( \frac{1}{%
n\ep^{d/2}} \right) \to q_2, \\
\Var(\Qt_n)& = \binom{n}{2}^{-2}
\bep^{-2}\sigma_n^2 = \Or \left(\frac{1}{n^2\ep^d}
\right) +\Or \left(\frac{1}{n} \right) +\oor\left(\frac{1}{n^3
\ep^{3d/2}}\right) +\oor\left(\frac{1}{n^2 \ep^{d}}\right) \notag \\
&=  \Or \left(\frac{1}{n^2\ep^d}
\right) +\Or \left(\frac{1}{n} \right)  \to 0 \ninf.  \notag
\end{align}
Hence, the assertion follows.\newline
\newline
\noindent \textit{(ii)} In order to avoid condition $\mc{A}$, we repeat
the argument of Proposition \ref{prop:Nn1}\textit{(i)} with the convergence
rates in Lemma \ref{lemma}\textit{(ii)-(iii)} replaced by the weaker limits
\begin{equation}  \label{eq:wl}
\sup_{\{i_1,i_2 \} \neq \{j_1,j_2\}} P(d(X_{i_1},X_{j_1}) \leq \ep,
d(X_{i_2},X_{j_2}) \leq \ep) \leq \sup_{i \neq j} P(d(X_i,X_j)\leq \ep)
\to 0 \mbox{ as } \ep \to 0,
\end{equation}
which follow from the stationarity and $m$-dependence of $\{X_i\}$ and since $%
P(X_i=X_j) =0$. First we obtain from \eqref{eq:ENn} that
\begin{equation}  \label{eq:E}
\mu_n= \binom{n}{2}q_{2,\ep} + \oor(n) \ninf.
\end{equation}
Moreover, if we use \eqref{eq:wl} in place of Lemma \ref{lemma} in the derivation of %
\eqref{VarN*}, \eqref{VarN2}, \eqref{V(N*-N)}, and thus finally %
\eqref{VarNn}, it follows that
\begin{align}  \label{eq:V}
&\Var(N^*_{n}) = \frac{n^2}{2} q_{2,\ep} + n^3
\left(\sigma^2_{1,0,\ep} + 2\sum_{h=1}^m \sigma^2_{1,h,\ep}
\right)+ \oor(n^2), \\
& \Var(N_{n}) = \Var(N^*_{n}) - \Var(N_{n}- N^*_{n}) + \oor(n^2), \notag \\
&\Var(N_{n}-N^*_{n}) = \oor(n^2),  \notag \\
&\sigma_n^2 = \frac{n^2}{2} q_{2,\ep} + n^3
\left(\sigma^2_{1,0,\ep} + 2\sum_{h=1}^m \sigma^2_{1,h,\ep}
\right)+ \oor(n^2), \ninf.  \notag
\end{align}
From \eqref{eq:E}, the last statement in \eqref{eq:V}, and the
condition $n\ep^d \to a, 0<a \leq \infty$,
\begin{align*}
\Ex\Qt_n & = \qt_{2,\ep} + \oor\left( \frac{1}{n\ep^{d}} \right) \to q_2, \\
\Var(\Qt_n)& = \binom{n}{2}^{-2}
\bep^{-2}\sigma_n^2 = \Or \left(\frac{1}{n^2\ep^d}
\right) +\Or \left(\frac{1}{n} \right) +\oor\left(\frac{1}{n^2
\ep^{2d}}\right) \to 0 \ninf,
\end{align*}
so the claim holds true. This completes the proof. \hfill $\Box$ \newline
\newline

\noindent \noindent \textit{Proof of Theorem \ref{th:preLim}.} \textit{(i)} Let

\begin{equation}  \label{eq:preLim1}
\sqrt{n}(\Qt_n-\qt_{2,\ep}) = \sqrt{n} \binom{n}{2}^{-1}
\bep^{-1}(N_{n} - \Ex N_{n}) + R_n,
\end{equation}
where, by Proposition \ref{prop:Nn1}\textit{(i)} and the assumption $n
\ep^d \to a, 0< a \leq \infty$,
\begin{equation}  \label{eq:preLim2}
R_n:=\sqrt{n} \binom{n}{2}^{-1} \bep^{-1}(\Ex N_{n} -
\binom{n}{2} q_{2,\ep}) = \oor \left(\frac{1}{\sqrt{n\ep^d}}
\right) \to 0 \ninf.
\end{equation}
Furthermore, from Proposition \ref{prop:Nn1}\textit{(ii)},
\begin{equation}  \label{eq:preLim3}
\left(\sqrt{n} \binom{n}{2}^{-1} \bep^{-1} \right)^2 \sigma_n^2
\to \nu/a+4\zeta_{1,m} \ninf.
\end{equation}
Finally, combining \eqref{eq:preLim1}, \eqref{eq:preLim2}, \eqref{eq:preLim3}%
, Theorem \ref{th:Nn2}\textit{(iii)}, and the Slutsky theorem gives the assertion
for $\Qt_n$. The statement about $\tilde{H}_n$ follows from Proposition 2 in Leonenko and Seleznjev (2010). \\

\noindent\textit{(ii)} The details are omitted, since the argument is similar to that of {\it (i)}, using the decomposition
corresponding to \eqref{eq:preLim1} with the $n\ep^{d/2}$-scaling. This completes the proof. \hfill $\Box$ \newline
\newline

\noindent \textit{Proof of Theorem \ref{th:Cons}.} {\it (i)}
As in Leonenko and Seleznjev (2010), the density smoothness condition yields
\begin{equation}\label{BiasBound}
|\qt_{2,\ep} -q_2| \leq
\frac{1}{2}K^2\ep^{2\alpha}.
\end{equation}
This bound and Proposition \ref{prop:Nn1}\textit{(i)} imply the assertion
\begin{equation*}
|\Ex \Qt_n - q_2| \leq |\qt_{2,\ep} -q_2| + |\mathrm{%
E} \Qt_n - \qt_{2,\ep}| \leq \frac{1}{2}K^2\ep^{2\alpha} + \mathrm{o%
}(1/(n\ep^{d/2})) \ninf.
\end{equation*}

\noindent \textit{(ii)}
Note that, by the assumptions $\ep \sim cn^{-2/(4\alpha +d)}$ and $0<\alpha \leq d/4$,
\begin{equation}\label{n2ep}
n^2\ep^d \sim c^dn^{8\alpha/(4\alpha+d)} \leq c^d n \ninf
\end{equation}
and hence, from \eqref{eq:VarQtn},
\begin{equation*}
\Var(\Qt_n) = \Or (n^{-8\alpha/(4\alpha+d)} )
\mbox{ as
} n \to \infty.
\end{equation*}
Further, since $\ep^{2\alpha} \sim c^{2\alpha}n^{-4\alpha/(4\alpha+d)}$, we get from \textit{(i)} and \eqref{n2ep} that the bias fulfills
\begin{align*}
|\Ex\Qt_n-q_2|
=\Or(n^{-4\alpha/(4\alpha+d)}) \ninf.
\end{align*}
Consequently, for some $C >0$ and any $A >0$,
\begin{equation*}
P(|\Qt_n-q_2| > An^{-4\alpha/(4\alpha+d)}) \leq n^{8\alpha/(4\alpha+d)} \frac{\Var(%
\Qt_n) + (\Ex \Qt_n - q_2)^2}{A^2 }
\leq \frac{C}{A^2}
\end{equation*}
and the desired convergence for $\Qt_n$ follows.
Moreover, combining this with Proposition 2 in Leonenko and Seleznjev (2010) proves the statement for $\tilde{H}_n$.\\

\noindent \textit{(iii)} The argument is similar to that  of {\it (ii)} and therefore is left out.
This completes the proof. \hfill $\Box$ \newline
\newline

\noindent \textit{Proof of Proposition \ref{std}.} First we study
the expectation of $U_{h,n}$. By Lemma \ref{lemma:lee}, the number of
3-tuples $(s_1,s_2,s_3)$ that satisfy $s_{i+1}-s_{i} > h+m, i = 1,2 $, is $%
\binom{n-2(m+h)}{3} $. Furthermore, observe that $3!$ of the elements $%
(i,j,k) \in \mc{E}_{h,n}$ are permutations of $(s_1,s_2,s_3)$. For the
corresponding variables, $X_j$ and $X_k$ are mutually independent and also
independent of $\{X_i,X_{i+h}\}$, so for these we obtain
\begin{equation*}
\Ex I(d(X_i,X_j)\leq \ep_0,d(X_{i+h},X_k)\leq \ep_0) =
\Ex p_{X,\ep_0}(X_i)p_{X,\ep_0}(X_{i+h})=\Ex p_{X,\ep_0}(X_1)p_{X,\ep_0}(X_{1+h}).
\end{equation*}
Thus, from Lemma \ref{lemma} and the assumption $n\ep_0^d \to \infty$,
\begin{align}  \label{std:bias}
\Ex U_{h,n} & = 3!\binom{n-2(m+2h)}{3} M_{h,n}^{-1}
b_{\ep_0}(d)^{-2}\Ex p_{X,\ep_0}(X_1)p_{X,%
\ep_0}(X_{1+h}) + \oor(1/(n\ep_0^d)) \\
& \to \Ex p(X_1)p(X_{1+h}) = q_{3,h} \ninf.  \notag
\end{align}
Next we consider the variance of $U_{h,n}$. Using the notation $%
I_{\ep_0} (i,j):=I(d(X_i,X_j)\leq \ep_0)$, we have
\begin{align}  \label{std:var}
\Var(U_{h,n} ) = M_{h,n}^{-2}b_{\ep_0}(d)^{-4} \!\!\!\!\!\!
\sum _{\substack{ (s_1,s_2,s_3) \in \mc{E}_{h,n}  \\ (t_1,t_2,t_3) \in
\mc{E}_{h,n}}} \!\!\!\!\!\! \Cov \left(
I_{\ep_0}(s_1,s_2)I_{\ep_0}(s_1+h,s_3),
I_{\ep_0}(t_1,t_2)I_{\ep_0}(t_1+h,t_3) \right).
\end{align}
We count the number of terms in this sum that are zero. Lemma \ref%
{lemma:lee} implies that the number of 6-tuples $\{u_1,\ldots, u_6 \}
\subseteq \{1,\ldots,n \}$ with $u_{i+1} - u_{i} > m+h, i = 1, \ldots, 5$,
is $\binom{n-5(m+h)}{6}$. Each such 6-tuple can be divided and permuted into
$\binom{6}{3}\times3!\times3! = 6!$ pairs $(s_1,s_2,s_3),(t_1,t_2,t_3) \in
\mc{E}_{h,n}$. The $m$-dependence of $\{X_i\}$ yields that the
corresponding random variables are independent, and hence at least
\begin{equation*}
6!\binom{n-5(m+h)}{6} = M_{h,n}^2 + \Or(n^5) \ninf
\end{equation*}
summands in \eqref{std:var} are zero. For each of the $\Or(n^5)$
non-zero terms,
\begin{align*}
|\Cov  ( I_{\ep_0}(s_1,s_2)I_{\ep_0}(s_1 +h,s_3)&,
I_{\ep_0}(t_1,t_2)I_{\ep_0}(t_1+h,t_3) ) |  \\
& \leq (\Ex I_{\ep_0}(s_1,s_2)I_{\ep_0}(s_1 +h,s_3))^{1/2} (\Ex I_{\ep_0}(t_1,t_2)I_{\ep_0}(t_1+h,t_3))^{1/2},
\end{align*}
so Lemma \ref{lemma}\textit{(iii)} gives that the sum of the
non-zero terms in \eqref{std:var} is $\oor(n^5\ep_0^d) %
\ninf$. Combining this with the condition $%
n\ep_0^{3d} \to c ,0<c \leq \infty$, we get
\begin{align}  \label{std:varAsy}
\Var(U_{h,n} ) = \oor \left(\frac{1}{n\ep_0^{3d}} \right)
\to 0 \ninf.
\end{align}
Finally, from \eqref{std:bias} and \eqref{std:varAsy} it follows that $\Ex(U_{h,n}- q_{3,h})^2 \to 0$, which completes the proof. \hfill $%
\Box$\newline
\newline

\noindent \textit{Proof of Theorem \ref{th:Norm1}.} The argument is similar
to that of Theorem 6 in Leonenko and Seleznjev (2010), so we show the main steps
only. From the decomposition
\begin{equation}  \label{eq:dec}
\sqrt{n}(\Qt_n-q_2) = \sqrt{n}(\Qt_n- \qt_{2,\ep}) +
\sqrt{n}(\qt_{2,\ep}-q_2),
\end{equation}
we see that the assertion for $\Qt_n$ is implied by the Slutsky theorem if
$\sqrt{n}(\Qt_n- \qt_{2,\ep}) \overset{\mathrm{D}}{\to}
N(0,\nu/a + \zeta_{1,m})$ and $|\sqrt{n}(\qt_{2,\ep}-q_2)| \to 0$%
. The asymptotic normality follows straight away from Theorem \ref{th:preLim}%
\textit{(i)}. Furthermore, the conditions for $\alpha$ and $\ep$
together with bound \eqref{BiasBound} lead to the desired convergence of $|\sqrt{n}(\qt_{2,\ep}-q_2)|$.
Finally, Proposition 2 in Leonenko and Seleznjev (2010) proves the claim for $\tilde{H}_n$.
This completes the proof. \hfill $\Box$ \newline
\newline

\noindent \textit{Proof of Theorem \ref{th:Norm2}.} \noindent \textit{(i)}
We use the decomposition corresponding to \eqref{eq:dec}:
\begin{equation}  \label{eq:Norm2:dec}
n^{\beta/2} c^{d/2}(\Qt_{n}- q_2) = n^{\beta/2} c^{d/2}(\Qt_{n} -
\qt_{2,\ep}) + n^{\beta/2} c^{d/2}(\qt_{2,\ep}-q_2).
\end{equation}
Note that the condition $\ep \sim cn^{-(2-\beta)/d}$ gives $n\ep^d \to 0$ and $n^{\beta/2} c^{d/2} /(n\ep^{d/2}) \to 1$, so the asymptotic
normality
\begin{equation}\label{eq:Norm2:norm}
n^{\beta/2} c^{d/2}(\Qt_n-\qt_{2,\ep}) \overset{%
\mathrm{D}}{\to} N(0, \nu) \ninf
\end{equation}
follows from Theorem \ref{th:preLim}\textit{(ii)} and the Slutsky theorem. %
Further, the assumptions $\alpha > (d/4)C_\beta$, $\ep \sim cn^{-(2-\beta)/d}$, and bound \eqref{BiasBound} imply
\begin{equation}  \label{eq:Norm2:bias}
|n^{\beta/2} c^{d/2}(\qt_{2,\ep}-q_2)| \leq
c^{d/2}\frac{1}{2}K^2n^{\beta/2}\ep^{2\alpha} \sim c^{d/2+2\alpha}\frac{1}{2}K^2n^{\beta/2 - 2\alpha(2-\beta)/d} \to 0 \ninf,
\end{equation}
since $\beta/2 - 2\alpha(2-\beta)/d < \beta/2- 2(d/4)C_\beta(2-\beta)/d=0$.
Thus, from \eqref{eq:Norm2:dec}, \eqref{eq:Norm2:norm}, \eqref{eq:Norm2:bias}, and the Slutsky theorem,
we obtain the statement for $\Qt_n$.
The assertion for $\tilde{H}_n$ follows by an argument similar to that of Proposition 2 in Leonenko and Seleznjev (2010). \\\\
\noindent \textit{(ii)} The argument follows the same steps as that of {\it (i)} and consequently is omitted. This completes the proof. \hfill $\Box$
\newline

\newpage {\noindent {\large \textbf{References}}}

\begin{reflist}

{\small  Anderson, T.W.\ (1971), {\it The Statistical Analysis of Time Series}, New York: John Wiley and Sons.}

{\small Barbour, A.D., Holst, L., Janson, S.\ (1992), {\it Poisson Approximation}, Oxford: Oxford University Press.}

{\small Baryshnikov, Y., Penrose, M.D., Yukich, J.E.\ (2009), 'Gaussian limits for generalized spacings', {\it Ann.\ Appl.\ Probab.},
19, 158--185.}

{\small Bickel, P.J. and Ritov, Y.\ (1988), 'Estimating integrated squared density derivatives: sharp best order of convergence estimates',
{\it Sankhy{\=a}: The Indian Journal of Statistics}, Series A, 381--393.}

{\small Billlingsley, P.\ (1995), {\it Probability and Measure}, New York: Wiley.}

{\small Escolano, F., Suau, P., Bonev, B.\ (2009), {\it Information Theory in Computer Vision and Pattern Recognition}, New York: Springer.}

{\small Gin\'{e}, E., Nickl, R.\ (2008), 'A simple adaptive estimator for the integrated square of a density', {\it Bernoulli}, 14, 47--61.}

{\small Goria, M.N., Leonenko, N.N., Mergel, V.V., Inverardi, P.L.N.\ (2005),
'A new class of random vector entropy estimators and its applications in testing statistical hypotheses', {\it J.\ Nonparam.\ Stat.}, 17, 277--297.}

{\small Gregorio, A., Iacus, S.M.\ (2009), 'On R\'{e}nyi information for ergodic diffusion processes', {\it Inform.\ Sci.}, 179, 279--291.}

{\small Harrelson, D., Houdr\'{e}, C.\ (2003), 'A characterization of $m$-dependent stationary infinitely divisible sequences with applications to weak
convergence', {\it Ann.\ Probab.}, 31, 849--881.}

{\small Hart, J.D., Vieu, P.\ (1990), 'Data-driven bandwidth choice for density estimation based on dependent data',
{\it Ann.\ Stat.}, 18, 873--890.}

{\small Jammalamadaka, S.R., Janson, S.\ (1986), 'Limit theorems for a
triangular scheme of $U$-statistics with applications to inter-point
distances', {\it Ann.\ Probab.}, 14, 1347-1358.}

{\small Johnson, O., Vignat C.\ (2007), 'Some results concerning maximum R\'{e}nyi entropy distributions',
{\it Ann.\ Inst.\ H.\ Poincar\'{e} Probab.\ Statist.}, 43, 339--351.}

{\small Joe, H.\ (1997), {\it Multivariate Models and Dependence Concepts}, London: Chapman and Hall.}

{\small Kapur, J.N.\ (1989), {\it Maximum-entropy Models in Science and Engineering}, New York: Wiley.}

{\small Kapur, J.N., Kesavan, H.K.\ (1992), {\it Entropy Optimization Principles with Applications}, New York: Academic Press.}

{\small Kim, T.Y., Luo, Z.M., Kim, C.\ (2011), 'The central limit theorem for
degenerate variable $U$-statistics under dependence', {\it J.\ Nonparam.\ Stat.},
23, 683-699.}

{\small Koroljuk, V.S., \ Borovskich, Y.V.\ (1994), {\it Theory of $U$-statistics}, Dordrecht: Kluwer.}

{\small Kotz, S., Balakrishnan, N., Johnson, N.L.\ (2000), {\it Continuous multivariate distributions: Models and applications},
New York: Wiley.}

{\small K\"{a}llberg, D., Seleznjev, O.\ (2012), 'Estimation of entropy-type integral functionals', preprint arXiv:1209.2544.}

{\small K\"{a}llberg, D., Leonenko, N., Seleznjev, O.\ (2012), 'Statistical inference for R\'{e}nyi entropy functionals',
{\it Lecture Notes in Comput.\ Sci.}, 7260, 36-51.}

{\small Last, Y., Klein, M., Kandel, A.\ (2001), 'Knowledge discovery in time
series databases', {\it IEEE Trans.\ Systems.\ Man.\ and Cybernetics - Part B}, 31, 160--169.}

{\small Laurent, B.\ (1996), 'Efficient estimation of integral functionals of a density',
{\it Ann.\ Statist.}, 24, 659--681.}

{\small Lee, A.J.\ (1990), {\it $U$-Statistics: Theory and Practice}, New York: Marcel Dekker.}

{\small Leonenko, N., Pronzato, L., Savani, V.\ (2008), 'A class of R\'{e}nyi
information estimators for multidimensional densities', {\it Ann.\ Statist.}, 36,
2153-2182. Corrections, (2010), {\it Ann.\ Statist.}, 38, 3837-3838.}

{\small Leonenko, N., Seleznjev, O.\ (2010), 'Statistical inference for the $%
\ep$-entropy and the quadratic R\'{e}nyi entropy', {\it J.\ Multivariate
Anal.}, 101, 1981-1994.}

{\small Neemuchwala, H., Hero, A., Carson, P.\ (2005), 'Image matching using
alpha-entropy measures and entropic graphs', {\it Signal Processing}, 85, 277-296.}

{\small Pardo, L.\ (2006), {\it Statistical Inference Based on Divergence Measures}, Boca Raton: Chapman \& Hall.}

{\small Penrose, M., Yukich, J.E.\ (2011), 'Limit theory for point processes in
manifolds', {\it Annals of Applied Probability}, to appear, see also preprint arXiv:1104.0914v1.}

{\small Principe, J.C.\ (2010), {\it Information Theoretic Learning}, New
York: Springer.}

{\small R\'{e}nyi, A.\ (1970), {\it Probability Theory}, Amsterdam: North-Holland.}

{\small Seleznjev, O., Thalheim, B.\ (2003), 'Average case analysis in database
problems', {\it Methodol.\ Comput.\ Appl.\ Prob.}, 5, 395-418.}

{\small Seleznjev, O., Thalheim, B.\ (2010), 'Random databases with approximate
record matching', {\it Methodol.\ Comput.\ Appl.\ Prob.}, 12, 63--89. }

{\small Serfling, R.J.\ (2002), {\it Approximation Theorems of Mathematical
Statistics}, New York: Wiley.}

{\small Shannon, C.E.\ (1948), 'A mathematical theory of communication', {\it Bell
Syst.\ Tech.\ Jour.}, 27, 379-423, 623--656.}

{\small Thalheim, B.\ (2000), {\it Entity-Relationship Modeling. Foundations of
Database Technology}, Berlin: Springer-Verlag.}

{\small Tsallis, C.\ (1988), 'Possible generalization of Boltzmann-Gibbs statistics',
{\it J.\ Stat.\ Phys.}, 52, 479--487.}

{\small Ullah, A.\ (1996), 'Entropy, divergence and distance measures with
econometric applications'. {\it J.\ Statist.\ Plann.\ Inference}, 49, 137--162.}

{\small Vatutin,V. A., Mikhailov, V.G.\ (1995), 'Statistical estimation of the entropy of discrete random variables with a large number of outcomes', {\it Russian Math.\ Surveys}, 50, 963-976.}

{\small Wang, Q.\ (1999), 'On Berry-Esseen rates for $m$-dependent $U$%
-statistics', {\it Stat.\ Prob.\ Letters}, 41, 123-130.}

\end{reflist}

\section{Appendix. Estimation for discrete distributions}

Consider a stationary $m$-dependent sequence $\{X_1,\ldots, X_n\}$ with
discrete $d$-dimensional (marginal) distribution $\mc{P}= \{p(k), k \in
N^d\}$. We present some results on the estimation of quadratic R\'{e}nyi
entropy for discrete distributions
\begin{equation*}
h_2(\mc{P}) := -\log \left( \sum_{k} p(k)^2 \right) = -\log(q_2)
\end{equation*}
and the corresponding quadratic functional
\begin{equation*}
q_2 := \sum_k p(k)^2 = P(X=Y),
\end{equation*}
where $X$ and $Y$ are independent vectors with distribution $\mc{P}$.
Similarly to the continuous case, let
\begin{equation*}
\zeta _{1,m}:=\Var(p(X_1))+2\sum_{h=1}^{m} \Cov%
(p(X_1),p(X_{h+t})).
\end{equation*}
Define the normalized statistic
\begin{equation*}
Q_{n}:=\binom{n}{2}^{-1}\sum_{1\leq i<j\leq n}I(X_{i}=X_{j}) = \binom{n}{2}^{-1} N_{n,0},
\end{equation*}%
to be an estimator for $q_{2}$. Let $H_{n}:=-\log (\max (Q_{n},1/n))$ be the corresponding estimator for $h_{2}$.
For $h=0,\ldots,m$, we also introduce the following estimator for $q_{3,h}:=\Ex p(X_1)p(X_{1+h})$,
\begin{equation*}
U_{h,n}:=M_{h,n}^{-1}\sum_{(i,j,k)\in \mc{E}%
_{h,n}}I(X_{i}=X_{j}=X_{i+h}=X_{k}),
\end{equation*}%
where $\mc{E}_{n,h}$ and $M_{n,h}$ are defined as in Section \ref%
{sec:Main}. By an argument similar to that of Proposition~\ref{std}, we get $%
U_{h,n}\overset{\mathrm{P}}{\rightarrow }q_{3,h}$. Hence, for $r\geq m$, a
consistent estimator for $\zeta _{1,m}$ is given by
\begin{equation*}
s_{r,n}^{2}:=U_{0,n}-Q_{n}^{2}+2\sum_{i=1}^{r}(U_{h,n}-Q_{n}^{2}).
\end{equation*}
Some asymptotic properties for the estimators of $q_{2}$ and $h_{2}$ follow by
combining the results of Ch.~2 in Lee (1990), Wang (1999), and the Slutsky
theorem.

\begin{theorem} \

\begin{itemize}
\item[(i)] For the expectation and variance, we obtain
\begin{align*}
\Ex (Q_n) & = q_2 + \Or(n^{-1}), \\
\Var(Q_n) & = 4\zeta_{1,m}n^{-1} + \Or(n^{-2}) \mbox{ as } n
\to \infty,
\end{align*}
and thus $Q_n$ and $H_n$ are consistent estimators for $q_2$ and $h_2$,
respectively.

\item[(ii)] If $\zeta_{1,m} >0$ and $r \geq m$, then
\begin{align*}
&\sqrt{n}(Q_n-q_2) \overset{\mathrm{D}}{\to} N(0, 4\zeta_{1,m}) \quad%
\mbox{and}\quad \frac{\sqrt{n}}{2s_{r,n}}(Q_n - q_2) \overset{\mathrm{D}}{\to}
N(0,1); \\
& \frac{\sqrt{n}}{2s_{r,n}}Q_n(H_n - h_2) \overset{\mathrm{D}}{\to} N(0, 1) %
\ninf.
\end{align*}
\end{itemize}
\end{theorem}

As in the continuous case, we have $\zeta_{1,m} \geq 0$ with equality, e.g., if $\mc P$ is uniform.
\end{document}